\documentclass[11pt]{amsart}
\usepackage[margin=3.2cm]{geometry}
\usepackage{amssymb,amsthm}
\usepackage{enumitem}
\usepackage[all]{xy}

\usepackage{etex} 
\usepackage{pictex} 

\allowdisplaybreaks
\numberwithin{equation}{subsection}

\theoremstyle{plain}

\newtheorem*{lem*}{Lemma}
\newtheorem*{prop*}{Proposition}
\newtheorem*{thm*}{Theorem}
\newtheorem*{cor*}{Corollary}
\newtheorem*{conj*}{Conjecture}
\theoremstyle{remark}
\newtheorem{prob}{Problem}

\newtheorem*{rmk}{Remark}

\newtheorem*{example}{Example}
\newtheorem*{examples}{Examples}

\newcommand{\bil}[2]{\langle #1, #2 \rangle}

\newcommand{\Ext}{\operatorname{Ext}}
\newcommand{\Hom}{\operatorname{Hom}}

\newcommand{\ind}{\operatorname{ind}}
\renewcommand{\le}{\leqslant}
\renewcommand{\ge}{\geqslant}
\newcommand{\GL}{\mathrm{GL}}
\newcommand{\SL}{\mathrm{SL}}

\newcommand{\T}{\mathrm{T}}

\renewcommand{\L}{\mathrm{L}}
\newcommand{\M}{\mathrm{M}}

\newcommand{\F}{\mathfrak{F}}
\newcommand{\St}{\mathrm{St}}
\newcommand{\ch}{\operatorname{ch}}

\newcommand{\soc}{\operatorname{soc}}
\newcommand{\hatQ}{\widehat{\mathrm{Q}}}

\begin{document}
\title[Decomposition of tensor products]{Decomposition of tensor
products of modular irreducible representations for $\SL_3$\\
(With an Appendix by C.M.~Ringel)}
\author{C.~Bowman} 
\address{CB: Corpus Christi College, Cambridge, CB2 1RH, England, UK}
\email{C.Bowman@dpmms.cam.ac.uk}
\author{S.R.~Doty} 
\address{SRD: Mathematics and Statistics, Loyola
University Chicago, Chicago, IL 60626, USA}
\email{doty@math.luc.edu} 
\author{S.~Martin}
\address{SM: Magdalene College, Cambridge, CB3 0AG, England, UK}
\email{S.Martin@dpmms.cam.ac.uk}
\address{CMR: Fakult\"{a}t f\"{u}r Mathematik,
Universit\"{a}t Bielefeld,
PO Box 100 131,
D-33 501 Bielefeld, Germany} 
\email{ringel@math.uni-bielefeld.de}
\date{24 October 2010}
\begin{abstract}
We give an algorithm for working out the indecomposable direct
summands in a Krull--Schmidt decomposition of a tensor product of two
simple modules for $G=\SL_3$ in characteristics 2 and 3.  It is shown
that there is a finite family of modules such that every such
indecomposable summand is expressible as a twisted tensor product of
members of that family.

Along the way we obtain the submodule structure of various Weyl and
tilting modules. Some of the tilting modules that turn up in
characteristic 3 are not rigid; these seem to provide the first
example of non-rigid tilting modules for algebraic groups. These
non-rigid tilting modules lead to examples of non-rigid projective
indecomposable modules for Schur algebras, as shown in the Appendix.

Higher characteristics (for $\SL_3$) will be considered in a later
paper.
\end{abstract}

\maketitle

\section{Introduction}\noindent
We begin by explaining our motivation, which may be formulated for an
arbitrary semisimple algebraic group in positive characteristic.

\subsection{}\label{intro:gen}
Let $G$ be a semisimple, simply connected linear algebraic group over
an algebraically closed field $K$ of positive characteristic $p$. We
fix a Borel subgroup $B$ and a maximal torus $T$ with $T \subset B
\subset G$ and we let $B$ determine the \emph{negative} roots. We
write $X = X(T)$ for the character group of $T$ and let $X^+$ denote
the set of dominant weights. By $G$-module we always mean a rational
$G$-module, i.e.\ a $K[G]$-comodule, where $K[G]$ is the coordinate
algebra of $G$.  For each $\lambda \in X^+$ we have the following (see
\cite{Jantzen}) finite dimensional $G$-modules:
\[
\begin{tabular}{ll}
$\L(\lambda)$ & simple module of highest weight $\lambda$;\\
$\Delta(\lambda)$ & Weyl module of highest weight $\lambda$;\\
$\nabla(\lambda)$ & $= \ind_B^G K_\lambda$; dual Weyl module of
 highest weight $\lambda$; \\
$\T(\lambda)$ & indecomposable tilting module of highest weight
$\lambda$
\end{tabular}
\]
where $K_\lambda$ is the 1-dimensional $B$-module upon which $T$ acts
by the character $\lambda$ with the unipotent radical of $B$ acting
trivially.  The simple modules $\L(\lambda)$ are contravariantly
self-dual.  The module $\nabla(\lambda)$ has simple socle isomorphic
to $\L(\lambda)$; the module $\Delta(\lambda)$ is isomorphic to $^\tau
\nabla(\lambda)$, the contravariant dual of $\nabla(\lambda)$, hence
has simple head isomorphic to $\L(\lambda)$.  

The central problem which interests us is as follows.

\begin{prob}
Describe the indecomposable direct summands of an arbitrary tensor
product of the form $\L(\lambda) \otimes \L(\mu)$, for $\lambda, \mu
\in X^+$.
\end{prob}

As usual, a superscript $M^{[j]}$ on a $G$-module $M$ indicates that
the structure has been twisted by the $j$th power of the Frobenius
endomorphism on $G$.  By Steinberg's tensor product theorem, there are
twisted tensor product factorizations
\begin{align*}
 \L(\lambda) &\simeq \L(\lambda^0) \otimes \L(\lambda^1)^{[1]} \otimes
 \L(\lambda^2)^{[2]} \otimes \cdots; \\
 \L(\mu) &\simeq \L(\mu^0) \otimes \L(\mu^1)^{[1]} \otimes
 \L(\mu^2)^{[2]} \otimes \cdots
\end{align*}
where $\lambda = \sum \lambda^j p^j$, $\mu = \sum \mu^j p^j$ are the
$p$-adic expansions (unique) such that each $\lambda^j, \mu^j$ belongs
to the restricted region
\[
X_1 = \{ \nu \in X^+ \mid \langle \alpha^\vee, \nu \rangle \le
p-1 \text{ for all simple roots } \alpha \}.
\]
Putting these factorizations into the original tensor product we obtain
\begin{equation}\label{eq:stp}
 \L(\lambda) \otimes \L(\mu) \simeq \textstyle\bigotimes_{j \ge 0}\;
 \big(\L(\lambda^1) \otimes \L(\mu^1)\big)^{[j]}
\end{equation}
and thus we see that in Problem 1 one should first study the case
where both highest weights in question are restricted.

Assume that Problem 1 has been solved for all pairs of restricted
weights (note that this is a {\em finite} problem for any given $G$).
Let $\F = \F(G)$ be the set of isomorphism classes of indecomposable
direct summands appearing in some $\L \otimes \L'$, for a pair $L,
L'$ of restricted simple $G$-modules. Let $[L \otimes L': I]$ be the
multiplicity of $I \in \F$ as a direct summand of $L \otimes L'$. Then
one can express each tensor product $\L(\lambda^j) \otimes \L(\mu^j)$
as a finite direct sum of indecomposable modules
\begin{equation} \label{eq:A}
 \L(\lambda^j) \otimes \L(\mu^j) \simeq \bigoplus_{I \in \F}\;
   [\L(\lambda^j) \otimes \L(\mu^j) : I]\; I .
\end{equation}
Thus, the original tensor product $\L(\lambda) \otimes \L(\mu)$ has a
decomposition of the form
\begin{equation*}
\L(\lambda) \otimes \L(\mu) \simeq \textstyle \bigotimes_{j\ge 0} \;
\bigoplus_{I \in \F}\; [\L(\lambda^j) \otimes \L(\mu^j) : I]\; I^{[j]}
\end{equation*}
and by interchanging the order of the product and sum we obtain the
decomposition
\begin{equation} \label{eq:star}
\textstyle \L(\lambda) \otimes \L(\mu) \simeq \bigoplus\, \big(\prod_{j
  \ge 0}\, [\L(\lambda^j) \otimes \L(\mu^j) : I_{j}]\big)\;
  \bigotimes_{j\ge 0} I_{j}^{[j]} 
\end{equation}
where the direct sum is taken over the set of all finite sequences
$(I_{0}, I_{1}, I_{2}, \dots)$ of members of $\F$.

This gives a direct sum decomposition of $\L(\lambda) \otimes \L(\mu)$
in terms of twisted tensor products of modules in $\F$. If all such
twisted tensor products are themselves indecomposable as $G$-modules,
then we have in some sense solved Problem 1 for general $\lambda,
\mu$. Even when this isn't true we have still obtained a first
approximation towards a solution to Problem 1.  This leads us to the
following secondary set of problems:

\begin{prob}  
Given $G$,

  (a) classify the members of the family $\F = \F(G)$ and compute the
multiplicities $[L\otimes L': I]$ for $I \in \F$, $L,L'$ restricted;

  (b) determine conditions under which a twisted tensor product of
members from $\F$ remains indecomposable;

  (c) determine the module structure of the members of $\F$.
\end{prob}

\noindent
Let $G_r$ denote the kernel of the $r$th iterate of the Frobenius, let
$G_rT$ denote the inverse image of $T$ under the same map, and let
$\hatQ_r(\lambda)$ denote the $G_rT$-injective hull of $\L(\lambda)$
for any $\lambda \in X_r$, where
\[
  X_r := \{ \nu \in X^+ \mid \langle \alpha^\vee, \nu \rangle \le
p^r-1 \text{ for all simple roots } \alpha \}.
\]
Let $h$ denote the Coxeter number of $G$.  If $p \ge 2h-2$ then
$\hatQ_r(\mu)$ has (for any $\mu \in X_r$) a $G$-module structure; this
structure is unique in the sense that any two such $G$-module
structures are equivalent.  (These statements are expected to hold for
all $p$.)  Concerning Problem 2(b) we observe the following.

\begin{lem*}\label{lem:prob2c}
  Assume that $p \ge 2h-2$ or that if $p < 2h-2$ then $\hatQ_1(\mu)$
  has a unique $G$-module structure for all $\mu \in X_1$.  If each
  member of the sequence $(I_j)_{j\ge 0}$ ($I_j \in \F$) has simple
  $G_1T$-socle with restricted highest weight then the twisted tensor
  product $\bigotimes_{j \ge 0}\; I_{j}^{[j]}$ is indecomposable as a
  $G$-module. Hence $P$ is indecomposable.
\end{lem*}

\begin{proof}
  By assumption the socle of $I_j$ is simple, as a $G_1T$-module,
  hence has the form $L(\mu(j))$ for some $\mu(j) \in X_1$. Hence the
  module $I_j$ embeds in the $G_1T$-injective hull $\hatQ_1(\mu(j))$
  of $\L(\mu(j))$, for each $j$, so $P:= I_{0} \otimes I_{1}^{[1]}
  \otimes \cdots \otimes I_{m}^{[m]}$ embeds in $Q:= \hatQ_1(\mu(0))
  \otimes \hatQ_1(\mu(1))^{[1]} \otimes \cdots \otimes
  \hatQ_1(\mu(m))^{[m]}$.  By \cite[II.11.16 Remark 2]{Jantzen} the
  module $Q$ has a $G$-module structure and is isomorphic to
  $\hatQ_r(\mu)$, where $\mu = \sum_j\, \mu(t_j) p^j$. Since
  $\hatQ_r(\mu)$ has simple $G_rT$-socle $\L(\mu)$ it follows that $P$
  also has simple $G_rT$-socle $\L(\mu)$, and thus has simple
  $G$-socle $\L(\mu)$.
\end{proof}

We note that in Types $A_1$ and $A_2$ ($G = \SL_2$, $\SL_3$) it is
known that $\hatQ_1(\mu)$ has a unique $G$-module structure for all
$\mu \in X_1$, for any $p$.  In the case $G=\SL_2$ (studied in
\cite{DH}) it turns out that for any $p$ the members of $\F$ are
always indecomposable tilting modules with simple $G_1T$-socle of
restricted highest weight, so the determination of the family $\F$ and
the multiplicities $[L\otimes L': I]$ leads in that case to a complete
solution of Problem 1 for all pairs of dominant weights.  The purpose
of this paper is to examine the next most complicated case, namely the
case $G = \SL_3$. In that case, we will see that all members of $\F$
have simple $G_1T$-socle of restricted highest weight when $p=2$, and
this holds with only two exceptions when $p=3$, so the decomposition
\eqref{eq:star} is decisive in characteristic 2 and provides a great
deal of information in characteristic 3. 

Furthermore, although in characteristic 3 the summands in
\eqref{eq:star} are not always indecomposable, by analyzing the
further splittings which arise, we show that there is a finite
family $\F'$, closely related to $\F$, such that every indecomposable
direct summand of $\L(\lambda) \otimes \L(\mu)$ is isomorphic to a
twisted tensor product of members of $\F'$. Thus, we obtain a complete
solution to Problem 1 in characteristics 2 and 3.

\subsection{}
The paper is organized as follows.  In Section \ref{sec:gen} we recall
known facts that we use.  Our main technique is to compute structure
of certain Weyl modules (using a computer when necessary) and use that
structure to deduce structural information on certain tilting modules.
The main results obtained by our computations are given in Sections
\ref{sec:p2} and \ref{sec:p3}. To be specific, the structure of the
relevant Weyl modules is given in \ref{Weyl:p2} and \ref{Weyl:p3},
while the main results on tensor products --- including description of
the family $\F$, multiplicities $[L \otimes L': I]$ for restricted
simples $L,L'$ and $I \in \F$, and structure of members of $\F$ (in
most cases) --- are summarized in \ref{ss:rtpd} and \ref{main:3}. One
will also find worked examples in those sections.

In characteristic 2 all members of $\F(\SL_3)$ are tilting modules
with simple $G_1T$-socle of restricted highest weight, so the
decomposition \eqref{eq:star} gives a complete answer to Problem 1 for
all pairs of dominant weights.  This is similar to what happens for
$G=\SL_2$.  Moreover, each member of $\F(\SL_3)$ in this case is rigid
(a module is called \emph{rigid} if its radical and socle filtrations
coincide) and can be described by a strong diagram in the sense of
\cite{alp}. Recall that in \cite{alp} a module diagram is a directed
graph depicting the radical series of the module, in such a way that
vertices correspond to composition factors and edges to non-split
extensions, and a strong diagram is one in which the diagram also
determines the socle series. (One should consult \cite{alp} for
precise statements.)

Characteristic 3 is more complicated. (As standard notation, we write
$(a,b)$ for a highest weight of the form $a \varpi_1 + b\varpi_2$
where $\varpi_1, \varpi_2$ are the usual fundamental weights.)  First,
all but two of the members of $\F(\SL_3)$ have simple $G_1T$-socle of
restricted highest weight. The two exceptional cases are in fact
simple modules of highest weights $(5,2)$ and $(2,5)$ that are not
restricted, and so one is forced to consider possible further
splitting of summands in \eqref{eq:star}, in cases where one or both
of these modules appears in a twisted tensor product on the right hand
side. (This happens only if the tensor square of the Steinberg module
occurs in some factor in the right hand side of \eqref{eq:stp}.)  In
all cases those further splittings can be worked out; see Proposition
\ref{p3:decisive}. This leads to the finite family $\F'$ discussed in
the last paragraph of \ref{intro:gen}.

Furthermore, in characteristic 3 it turns out that four members of
$\F(\SL_3)$ --- namely the tilting modules $\T(3,3)$, $\T(4,3)$,
$\T(3,4)$, and $\T(4,4)$ --- are not rigid and do not have strong
Alperin diagrams.  The structure of one of the simplest of these
examples, $\T(4,3)$, is analyzed in detail in the Appendix by
C.M.~Ringel, using different methods. Although not itself projective,
Ringel shows that $\T(4,3)$ is a quotient of the corresponding
projective indecomposable for an appropriate Schur algebra, and thus
he produces an example of a non-rigid projective indecomposable module
for that Schur algebra. (See \cite{Green, Martin, Donkin:SA1,
  Donkin:SA2} for background on Schur algebras.) The other non-rigid
modules are subject to a similar analysis.

Preliminary calculations indicate that members of $\F(\SL_3)$ are
again rigid in characteristics higher than 3. The observed anomalies
in characteristic 3 are associated with the fact that some of the Weyl
modules which turn up are too close to the upper wall of the ``lowest
$p^2$-alcove'' and thus have composition factors with multiplicity
greater than 1. (Those multiplicities follow, e.g.\ from \cite{DS},
from knowledge of composition factor multiplicities in baby Verma
modules, which are well known in this case.) The simple characters for
$\SL_3$ have been known for a long time (see e.g., \cite{Jant1,
  Jant2}).

Our results overlap somewhat with \cite{jgj}, \cite{kh} although our
methods are different and we push the calculations further. Larger
characteristics, for which some calculations become in a sense
independent of $p$, will be treated in a future paper.

This paper has been circulating for some time in various forms, and
since the first version was made available, the preprint \cite{AK} has
appeared, in which further examples of non-rigid tilting modules for
algebraic groups are obtained.

\section{Preliminaries}\label{sec:gen}\noindent
We recall some general facts that will be used in our calculations.

\subsection{} \label{ss:genPillen}
Let us recall Pillen's Theorem \cite[\S2, Corollary A]{Pillen} (see
also \cite[Theorem (2.5)]{Donkin:Zeit}). Write $\St_r$ for the $r$th
Steinberg module $\L((p^r-1)\rho) = \Delta((p^r-1)\rho) =
\T((p^r-1)\rho)$. Then for $\lambda \in X_r$ the tilting module
$\T(2(p^r-1)\rho + w_0\lambda)$ is isomorphic to the indecomposable
$G$-component of $\St_r \otimes \L((p^r-1)\rho + w_0\lambda)$
containing the weight vectors of highest weight $2(p^r-1)\rho +
w_0\lambda$.

\subsection{}
In general the formal character of a tilting module is not known; even
for $\SL_3$, as far as we are aware this remains an open problem. The
following general result of Donkin (see \cite[Proposition
  5.5]{Donkin:TiltHandbook}) computes the formal character of certain
tilting modules. Let $\lambda, \mu \in X^+$ and assume that $(\lambda,
\alpha_0^\vee) \le p$, where $\alpha_0$ is the highest short
root. Then:
\begin{equation}\label{eq:Donkinch1}
  \ch \T((p-1) \rho + \lambda) = \ch \L((p-1)\rho) \sum_{\nu \in
    W\lambda} e(\nu)
\end{equation}
and for any $\nu \in X^+$,
\begin{equation}\label{eq:Donkinch2}
  (\T((p-1)\rho+\lambda+p\mu): \nabla(\nu)) = \sum_{\xi \in N(\nu)}
  (\T(\mu):\nabla(\xi))
\end{equation}
where $N(\nu) = \{\xi \in X^+: \nu+\rho -p(\xi+\rho) \in W\lambda \}$.
Furthermore, in Lemma 5 of Section 2.1 in \cite{Donkin:coho}, the
characters of the tilting modules which are projective and
indecomposable as $G_1$-modules are computed explicitly, for $G =
\SL_3$.

\subsection{}
Another useful general fact (that will be used repeatedly) is the
observation that tilting modules are contravariantly self-dual:
\begin{equation}\label{eq:self-duality}
  {}^{\tau} \T(\lambda) \simeq \T(\lambda)
\end{equation}
for all $\lambda \in X^+$.  This is because (by
\cite[II.2.13]{Jantzen}) contravariant duality interchanges
$\Delta(\mu)$ and $\nabla(\mu)$, so ${}^\tau\T(\lambda)$ is again
indecomposable tilting, of the same highest weight.

\subsection{} \label{gen:tp-thm}
Finally, there is a twisted tensor product theorem for tilting
modules, assuming that Donkin's conjecture \cite[Conjecture
  (2.2)]{Donkin:Zeit} is valid or that $p \ge 2h-2$. (It is well known
\cite[II.11.16, Remark 2]{Jantzen} that the conjecture is valid for
all $p$ in case $G = \SL_3$.) For our purposes, it is convenient to
reformulate the tensor product theorem in the following form. First we
observe that, given $\lambda \in X^+$ satisfying the condition
\begin{equation}\label{eqn:condition}
  \langle \lambda, \alpha^\vee \rangle \ge p-1, \quad \text{for all
    simple roots } \alpha, 
\end{equation}
there exist unique weights $\lambda'$, $\mu$ such that
\begin{equation}\label{eqn:p-adic-start}
  \lambda = \lambda' + p \mu, \qquad \lambda' \in (p-1)\rho+ X_1,\ \mu
  \in X^+ .
\end{equation}
This is easy to see: for each $\lambda_i$ in $\lambda = \sum \lambda_i
\varpi_i$ where the $\varpi_i$ are the fundamental weights, express
$\lambda_i -(p-1)$ (uniquely) in the form $\lambda_i-(p-1) = r_i + p
s_i$ with $0 \le r_i \le p-1$. Then set $\lambda' = (p-1)\rho + \sum
r_i \varpi_i$ and $\mu = \sum s_i \varpi_i$.

Now by induction on $m$ using \eqref{eqn:condition} and
\eqref{eqn:p-adic-start} one shows that every $\lambda \in X^+$ has a
unique expression in the form
\begin{equation} \label{eqn:p-adic}
  \lambda = \textstyle\sum_{j=0}^m a_j(\lambda) \, p^j
\end{equation}
with $a_0(\lambda), \dots, a_{m-1}(\lambda) \in (p-1)\rho + X_1$ and
$\langle a_m(\lambda), \alpha^\vee \rangle < p-1$ for at least one
simple root $\alpha$.

Given $\lambda \in X^+$, express $\lambda$ in the form
\eqref{eqn:p-adic}. Assume Donkin's conjecture holds if $p <
2h-2$. Then there is an isomorphism of $G$-modules
\begin{equation}\label{eq:tptilt}
  \T(\lambda) \simeq \bigotimes_{j=0}^m \T(a_j(\lambda))^{[j]}.
\end{equation}
To prove this one uses induction and \cite[Lemma II.E.9]{Jantzen}
(which is a slight reformulation of \cite[Proposition
  (2.1)]{Donkin:Zeit}).

\section{Results for $p=2$}\label{sec:p2}\noindent
For the rest of the paper we take $G = \SL_3$.  Conventions: Dominant
weights are written as ordered pairs $(a,b)$ of non-negative integers;
one should read $(a,b)$ as an abbreviation for $a \varpi_1 + b
\varpi_2$ where $\varpi_1, \varpi_2$ are the fundamental weights,
defined by the condition $\bil{\varpi_i}{\alpha_j^\vee} =
\delta_{ij}$.  \emph{When describing module structure, we shall always
  identify a simple module $\L(\lambda)$ with its highest weight
  $\lambda$.}  Whenever possible we will depict the structure by
giving an Alperin diagram (see \cite{alp} for definitions) with edges
directed downwards, except in the uniserial case, where we will write
$M = [L_s, L_{s-1}, \dots, L_1]$ for a module $M$ with unique
composition series $0 = M_0 \subset M_{1} \subset \cdots \subset
M_{s-1} \subset M_s = M$ such that $L_j \simeq M_j/M_{j-1}$ is simple
for each $j$.

\subsection{Structure of certain Weyl modules for $p=2$} \label{Weyl:p2}
The results given below were computer generated, using GAP \cite{GAP}
code available on the second author's web page. (Some cases are
obtainable from \cite{D:thesis}.)

The restricted region $X_1$ in this case consists of the weights of
the form $(a,b)$ with $0 \le a,b \le 1$, and we have 
\begin{gather*}
\Delta(0,0) = \L(0,0),\ \Delta(1,0) = \L(1,0), \\ \Delta(0,1) =
\L(0,1),\ \Delta(1,1) = \L(1,1).
\end{gather*}
These are all tilting modules. Thus it follows immediately that all
the members of $\F$ are tilting.

The structure of the other Weyl modules we need is depicted below. The
uniserial modules have structure
\begin{gather*}
 \Delta(2,0) = [(2,0), (0,1)],\ \Delta(0,2)=  [(0,2), (1,0)]\\
 \Delta(3,0) =  [(3,0), (0,0)],\  \Delta(0,3) = [(0,3), (0,0)] \\
 \Delta(2,1) = [(2,1), (0,2), (1,0)],\  \Delta(1,2) = [(1,2), (2,0), (0,1)]. 
\end{gather*}
Finally, the structure of $\Delta(2,2)$ is given by the diagram
\[
\Delta(2,2) = \framebox[34mm]{
\begin{minipage}{34mm}
\def\objectstyle{\scriptstyle}
\xymatrix@=6pt{
  & (2,2) \ar@{-}[dl] \ar@{-}[dr] & \\
(0,3) \ar@{-}[dr] & & (3,0) \ar@{-}[dl] \\
  & (0,0) &
}\end{minipage}
} \; .
\]
We worked these out using explicit calculations in the hyperalgebra,
by methods similar to those of \cite{ron, xi}.

\subsection{Restricted tensor product decompositions for $p=2$}
\label{ss:rtpd}
The indecomposable decompositions of restricted tensor products for
$p=2$ is as follows. (We omit any decomposition of the form
$\L(\lambda) \otimes \L(\mu)$ where one of $\lambda, \mu$ is zero.)
There is an involution on $G$-modules which on weights is the map
$\lambda \to -w_0(\lambda)$, where $w_0$ is the longest element of the
Weyl group. (In Type $A$ this comes from a graph automorphism of the
Dynkin diagram.) We refer to this involution as \emph{symmetry}, and
we will often omit calculations and results that can be obtained by
symmetry from a calculation or result already given.

\begin{prop*}
Suppose $p=2$. 
\begin{enumerate}[label={\rm(\alph*)},leftmargin=*,itemsep=0.5em]

\item The indecomposable direct summands of tensor products
of non-trivial restricted simple $\SL_3$-modules are as follows: {\rm

(1)\quad $\L(1,0)\otimes \L(1,0) \simeq \T(2,0)$;\quad $\L(0,1)\otimes
\L(0,1) \simeq \T(0,2)$;

(2)\quad $\L(1,0)\otimes \L(0,1) \simeq \T(1,1) \oplus \T(0,0)$; 

(3)\quad $\L(1,0)\otimes \L(1,1) \simeq \T(2,1)$;\quad $\L(0,1)\otimes
\L(1,1) \simeq \T(1,2)$;

(4)\quad $\L(1,1)\otimes \L(1,1) \simeq \T(2,2) \oplus 2 \T(1,1)$.}

\noindent
Thus the family $\F(\SL_3)$ is in this case given by $\F = \{ \T(a,b):
0 \le a,b \le 2\}$.

\item The structure of the uniserial members of $\F$ is given as
  follows:

$\T(0,0) = [(0,0)]$;\quad $\T(1,0) = [(1,0)]$;\quad $\T(1,1) =
  [(1,1)]$;

$\T(2,0) = [(0,1), (2,0), (0,1)]$. 

\noindent
The structure diagrams of $\T(2,1)$, $\T(2,2)$ are displayed below:

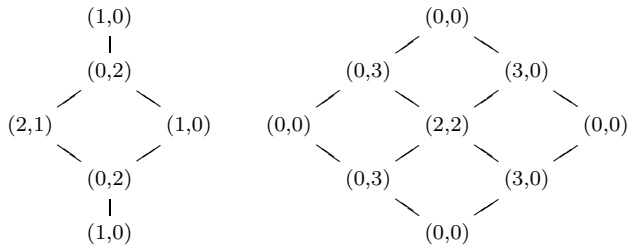
\begin{center}
\begin{tabular}{cc}
\def\objectstyle{\scriptstyle}
\xymatrix@=6pt{
  & (1,0) \ar@{-}[d] & \\
  & (0,2) \ar@{-}[dl] \ar@{-}[dr] & \\
(2,1) \ar@{-}[dr] & & (1,0) \ar@{-}[dl] \\
  & (0,2) \ar@{-}[d] & \\
  & (1,0) & 
}
&
\def\objectstyle{\scriptstyle}
\xymatrix@=6pt{
  & & (0,0) \ar@{-}[dl] \ar@{-}[dr] & & \\
  & (0,3) \ar@{-}[dl] \ar@{-}[dr] & & (3,0) \ar@{-}[dl] \ar@{-}[dr] & \\
(0,0) \ar@{-}[dr] & & (2,2) \ar@{-}[dl] \ar@{-}[dr] & & (0,0) \ar@{-}[dl] \\
  & (0,3) \ar@{-}[dr] & & (3,0) \ar@{-}[dl] & \\
  & & (0,0) & &
}
\end{tabular}
\end{center}
and the structure diagrams of $\T(0,1)$, $\T(0,2)$, and $\T(1,2)$ are
obtained by symmetry from cases already listed.

\item Each member of $\F$ has simple $G_1T$-socle (and head) with highest
  weight belonging to the restricted region $X_1$.
\end{enumerate}
\end{prop*}

\noindent
The proof is given in \ref{ss:pf}.  First we consider consequences and
give some examples. Recall that a dominant weight is called minuscule
if the weights of the corresponding Weyl module form a single Weyl
group orbit. For $G=\SL_3$ the minuscule weights are $(0,0)$, $(1,0)$,
and $(0,1)$.

\begin{cor*}
  Let $p=2$. Given arbitrary dominant weights $\lambda, \mu$ write
  $\lambda = \sum \lambda^j p^j$, $\mu = \sum \mu^j p^j$ with
  $\lambda^j, \mu^j \in X_1$ for all $j\ge 0$.
  \begin{enumerate}[label={\rm(\alph*)},leftmargin=*,itemsep=0.5em]

  \item In the decomposition \eqref{eq:star}, each term in the direct
    sum is indecomposable. Hence the indecomposable direct summands of
    $\L(\lambda) \otimes \L(\mu)$ are expressible as a twisted tensor
    product of members of $\F$. Conversely, every twisted tensor
    product of members of $\F$ occurs in some $\L(\lambda) \otimes
    \L(\mu)$.

  \item $\L(\lambda) \otimes \L(\mu)$ is indecomposable if and only if
    for each $j \ge 0$ the unordered pair $\{\lambda^j, \mu^j\}$ is
    one of the cases $\{(1,0),(1,0)\}$, $\{(0,1),(0,1)\}$,
    $\{(1,0),(1,1)\}$, $\{(0,1),(1,1)\}$ or one of $\lambda^j, \mu^j$
    is the trivial weight $(0,0)$.

  \item Let $m$ be the maximum $j$ such that at least one of
    $\lambda^j$, $\mu^j$ is non-zero. Then $\L(\lambda) \otimes
    \L(\mu)$ is indecomposable tilting, isomorphic to
    $\T(\lambda+\mu)$, if and only if: {\rm(i)} for each $0\le j \le
    m-1$, one of $\lambda^j, \mu^j$ is minuscule and the other is the
    Steinberg weight $(1,1)$, and {\rm(ii)} $\{\lambda^m, \mu^m\}$ is
    one of the cases listed in part {\rm(b)}.
  \end{enumerate}
\end{cor*}

\begin{proof}
  Part (a) follows from \eqref{eq:star} and Lemma
  \ref{intro:gen}. Part (b) follow from the proposition and the
  discussion preceding \eqref{eq:star}, which shows that each
  $\L(\lambda^j) \otimes \L(\mu^j)$ must be itself indecomposable in
  order for $\L(\lambda) \otimes \L(\mu)$ to be indecomposable. Then
  we get part (c) from part (b) by applying Donkin's tensor product
  theorem \eqref{eq:tptilt}.
\end{proof}

\begin{examples}
(i) To illustrate the procedure in part (a) of the corollary, we work
  out a specific example:
\begin{align*}
  \L(7,2) &\otimes \L(6,3)\\ &\simeq \big(\L(1,0)\otimes \L(0,1)\big)
  \otimes \big(\L(1,1)\otimes \L(1,1)\big)^{[1]} \otimes \big(\L(1,0)
  \otimes \L(1,0)\big)^{[2]} \\ &\simeq \big(\T(1,1)\oplus
  \T(0,0)\big) \otimes \big(\T(2,2)\oplus 2\T(1,1)\big)^{[1]} \otimes
  \T(2,0)^{[2]} \\ &\simeq \T(13,5) \oplus 2\T(6,2)^{[1]} \oplus
  2\T(11,3) \oplus 2\T(5,1)^{[1]} .
\end{align*}
In the calculation, the first line follows from Steinberg's tensor
product theorem, the second is from the proposition, and to get the
last line we applied Donkin's tensor product theorem
\eqref{eq:tptilt}, after interchanging the order of sums and products.

\noindent
(ii) We have $\L(3,0) \otimes \L(3,2) \simeq \big(\L(1,0) \otimes
\L(1,0)\big) \otimes \big(\L(1,0)\otimes \L(1,1)\big)^{[1]} \simeq
\T(2,0) \otimes \T(2,1)^{[1]}$, which is indecomposable but not
tilting. This illustrates the procedure in part (b) of the corollary.

\noindent
(iii) We have $\L(3,0)\otimes \L(3,1) \simeq \big(\L(1,0) \otimes
\L(1,1)\big) \otimes \big(\L(1,0)\otimes \L(1,0)\big)^{[1]} \simeq
\T(2,1) \otimes \T(2,0)^{[1]} \simeq \T(6,1)$, illustrating part (c)
of the corollary.

\noindent
(iv) It is not the case that every indecomposable tilting module occurs
as a direct summand of some tensor product of two simple modules. For
instance, neither $\T(3,0)$ nor $\T(0,3)$ (both of which are uniserial
of length 3) can appear as one of the indecomposable direct summands on
the right hand side of \eqref{eq:star}. This follows from
\eqref{eq:tptilt}. More generally, this applies to any non-simple
tilting module of the form $\T(a,b)$ with one of $a,b$ equal to zero
and the other greater than 2.
\end{examples}

\subsection{}\label{ss:pf}
We now consider the proof of Proposition \ref{ss:rtpd}.  First we
compute the composition factor multiplicities of the restricted tensor
products. Let $\chi_p(\lambda)$ be the formal character of
$\L(\lambda)$. Then:
\begin{enumerate}
\makeatletter
\renewcommand{\labelenumi}{(\theenumi) }
\makeatother
\item $\chi_p(1, 0) \cdot \chi_p(1, 0) =  \chi_p(2, 0) + 2 \chi_p(0, 1)$;

\item $\chi_p(1, 0) \cdot \chi_p(0, 1) =  \chi_p(1, 1) + \chi_p(0, 0)$;

\item $\chi_p(1, 0) \cdot \chi_p(1, 1) = \chi_p(2, 1) + 2\chi_p(0, 2) + 3
\chi_p(1, 0)$;

\item $\chi_p(0, 1) \cdot \chi_p(0, 1) =  \chi_p(0, 2) + 2 \chi_p(1, 0)$;

\item $\chi_p(0, 1) \cdot \chi_p(1, 1) = \chi_p(1, 2) + 2\chi_p
(2, 0) + 3\chi_p(0, 1)$;

\item $\chi_p(1, 1) \cdot \chi_p(1, 1) = \chi_p(2, 2) + 2 \chi_p(0, 3) +
2\chi_p(3, 0) + 2\chi_p(1, 1) + 4\chi_p(0, 0)$.
\end{enumerate}

Since $\L(1,0) = \T(1,0)$, it follows that $\L(1,0)\otimes \L(1,0)$ is
tilting. It must have $\T(2,0)$ as a direct summand by highest weight
considerations. But $\T(2,0)$ is contravariantly self-dual with
$\L(0,1)$ in the socle, so it follows that $\L(0,1)$ appears with
multiplicity at least 2 as a composition factor of $\T(2,0)$. Now
character considerations force the structure to be given by
\[
  \L(1,0)\otimes \L(1,0) \simeq \T(2,0)
\]
where $\T(2,0) = [(0,1), (2,0), (0,1)]$.  By symmetry we also have
\[
  \L(0,1)\otimes \L(0,1) \simeq \T(0,2)
\]
where $\T(0,2) = [(1,0), (0,2), (1,0)]$.  

$\L(1,0)\otimes \L(0,1)$ is tilting and has a direct summand
isomorphic to $\T(1,1) = \L(1,1)$.  By character considerations
it follows that there is one other indecomposable summand, namely
$\T(0,0) = \L(0,0)$. Hence
\[
\L(1,0)\otimes \L(0,1) \simeq \T(0,0) \oplus \T(1,1).
\]

$\L(1,0)\otimes \L(1,1)$ is tilting and has a direct summand
$\T(2,1)$. Self-duality of $\T(2,1)$ forces a copy of $\L(1,0)$ at the
top, extending $\L(0,2)$. This, along with the structure of the Weyl
modules and known Ext information forces the structure of $\T(2,1)$ to
be as given in the statement of Proposition \ref{ss:rtpd}(b), and also
forces
\[
\L(1,0)\otimes \L(1,1) \simeq \T(2,1).
\]
By symmetry we obtain also
\[
 \L(0,1)\otimes \L(1,1) \simeq \T(1,2). 
\]

Finally, $\L(1,1)\otimes \L(1,1)$ is tilting, with a direct summand
isomorphic to $\T(2,2)$.  The highest weights of all simple
composition factors of the tensor product are in the same linkage
class, excepting $(1,1)$, which appears with multiplicity 2. So two
copies of $\T(1,1)$ split off. Moreover, $\T(2,2)$ has a
submodule isomorphic to $\Delta(2,2)$, thus contains $\L(0,0)$ in the
socle. This forces another copy of $\L(0,0)$ at the top of $\T(2,2)$,
and this along with known Ext information and the structure of the
Weyl modules forces the structure of $\T(2,2)$ to be as given in
Proposition \ref{ss:rtpd}(b), and also
forces 
\[
\L(1,1)\otimes \L(1,1) \simeq \T(2,2) \oplus 2 \T(1,1).
\]

All the claims in Proposition \ref{ss:rtpd}(a), (b) are now clear.  It
remains to verify the claim in (c). It is known that Donkin's
conjecture holds for $\SL_3$, as discussed at the beginning of
\ref{gen:tp-thm}, so $\T((p-1)\rho+\lambda)$ is as a $G_1T$-module
isomorphic to $\hatQ_1((p-1)\rho + w_0 \lambda)$ for any $\lambda \in
X_1$.  Thus $\T(2,1)$, $\T(1,2)$, and $\T(2,2)$ each has a simple
$G_1T$-socle of restricted highest weight. For $\T(2,0)$ and $\T(0,2)$
one can argue by contradiction, using the fact \cite[Proposition
  (1.5)]{Donkin:Zeit} that truncation to an appropriate Levi subgroup
$L$ maps indecomposable tilting modules for $G$ onto indecomposable
tilting modules for $L$. Thus $\T(2,0)$ and $\T(0,2)$ truncate to
$\T(2)$ for $L\simeq \SL_2$, which is known to have simple
$L_1T$-socle and length three. If $\T(2,0)$ or $\T(0,2)$ did not have
simple $G_1T$-socle then the same would be true of the truncation,
since no composition factors are killed under truncation. Claim (c)
for the remaining cases is trivial.

\section{Results for $p=3$}\label{sec:p3}\noindent
In characteristic 3 several of the Weyl modules one must consider are
non-generic due to the proximity of their highest weight to the upper
wall of the lowest $p^2$-alcove. This leads ultimately to examples of
non-rigid tilting modules. Another complication is that the
$G_1T$-socles of two direct summands of the tensor square of the
Steinberg module fail to be simple.

\subsection{Structure of certain Weyl modules for $p=3$}\label{Weyl:p3}
We record the structure of certain Weyl modules needed later.  The
uniserial Weyl modules that turn up in our tensor product
decompositions have structure given by
\begin{gather*}
 \Delta(0,0) = \L(0,0),\ \Delta(1,0) = \L(1,0),\ 
 \Delta(2,0) = \L(2,0),\\   \Delta(2,1) = \L(2,1),\ 
 \Delta(2,2) = \L(2,2),\ \Delta(5,2) = \L(5,2),\\
 \Delta(1,1) = [(1,1), (0,0)],\  \Delta(3,0) = [(3,0), (1,1)],\\
 \Delta(4,0) = [(4,0), (0,2)],\ \Delta(3,1) = [(3,1), (1,2)], \\
 \Delta(5,0) = [(5,0), (0,1)],\ \Delta(5,1) = [(5,1), (1,0)], \\
 \Delta(3,2) = [(3,2), (1,3), (2,1)],\ \Delta(6,0) = [(6,0), (4,1), (0,0)]\\
 \Delta(4,2) = [(4,2), (0,4), (2,0)],\ \Delta(6,2) = 
[(6,2), (4,3), (1,0), (5,1)].
\end{gather*}
We note that the structure of $\Delta(6,2)$ is needed only in the
Appendix. The non-uniserial cases we need have structure
\[
  \Delta(4,1) = \framebox[34mm]{
\begin{minipage}{34mm}
\def\objectstyle{\scriptstyle}
\xymatrix@=6pt{
  & (4,1) \ar@{-}[dl] \ar@{-}[dr] \ar@{-}[d] & \\
(0,3) \ar@{-}[dr] & (0,0)\ar@{-}[d] & (3,0) \ar@{-}[dl] \\
  & (1,1) &
}\end{minipage}
}\, ,\qquad
  \Delta(4,3) = \framebox[34mm]{
\begin{minipage}{34mm}
\def\objectstyle{\scriptstyle}
\xymatrix@=6pt{
  & (4,3) \ar@{-}[d] & \\
  & (1,0) \ar@{-}[dl] \ar@{-}[dr]  & \\
(0,5) \ar@{-}[dr] &  & (5,1) \ar@{-}[dl]\\
 & (1,0) &
}\end{minipage}
} \, ,
\]
\[
  \Delta(3,3) = \framebox[34mm]{
\begin{minipage}{34mm}
\def\objectstyle{\scriptstyle}
\xymatrix@=6pt{
  & (3,3) \ar@{-}[d] & \\
  & (0,0) \ar@{-}[dl] \ar@{-}[dr]  & \\
(1,4) \ar@{-}[dr]\ar@{-}[d]\ar@{-}[drr] &  & (4,1) \ar@{-}[dl]\ar@{-}[d] 
 \ar@{-}[dll]\\
 (0,3)\ar@{-}[dr]  & (0,0)\ar@{-}[d] & (3,0)\ar@{-}[dl] \\
 & (1,1) &
}\end{minipage}
}\, ,\qquad
  \Delta(4,4) = \framebox[54mm]{
\begin{minipage}{34mm}
\def\objectstyle{\scriptstyle}
\xymatrix@=6pt{
 & & (4,4) \ar@{-}[ddl]\ar@{-}@/^1pc/[dddrr]
\ar@{-}@/_1pc/[dddll]\ar@{-}[ddr] & & \\
 \\
 & (3,3)\ar@{-}[dr] & & (1,1)\ar@{-}[dl]\ar@{-}[dll]\ar@{-}[d] & \\
 (0,6)\ar@{-}[dr] & (0,3)\ar@{-}[d]\ar@{-}[drr] & (0,0)\ar@{-}[dl]\ar@{-}[dr] & 
 (3,0)\ar@{-}[dll]\ar@{-}[d] & (6,0)\ar@{-}[dl] \\
 & (1,4)\ar@{-}[dr] & & (4,1)\ar@{-}[dl] & \\
 & & (0,0) & &
}\end{minipage}
} \,.
\]
As for the case $p=2$, these structures were obtained by explicit
calculations in the hyperalgebra, using GAP to do the calculations.

\subsection{Restricted tensor product decompositions for $p=3$} 
\label{main:3}
The indecomposable decompositions of restricted tensor products for
$p=3$ is given below. We omit any decomposition of the form
$\L(\lambda) \otimes \L(\mu)$ where one of $\lambda, \mu$ is zero, and
we omit all cases that follow by applying symmetry to a case already
listed.

\begin{prop*}
  Let $p=3$. 
\begin{enumerate}[label={\rm(\alph*)},leftmargin=*,itemsep=0.5em]
\item The indecomposable direct summands of tensor products of
non-trivial restricted simple $\SL_3$-modules are as follows:{\rm

(1)\quad $\L(1,0)\otimes \L(1,0) \simeq \T(2,0)\oplus \T(0,1)$;

(2)\quad $\L(1,0)\otimes \L(0,1) \simeq \T(1,1)$;

(3)\quad $\L(1,0)\otimes \L(2,0) \simeq \T(3,0)$;

(4)\quad $\L(1,0)\otimes \L(1,1) \simeq \T(2,1) \oplus \T(0,2)$;

(5)\quad $\L(1,0)\otimes \L(0,2) \simeq \T(1,2) \oplus \T(0,1)$;

(6)\quad $\L(1,0)\otimes \L(2,1) \simeq \T(3,1) \oplus  \T(2,0)$;

(7)\quad $\L(1,0) \otimes \L(1,2) \simeq \T(2,2) \oplus \T(0,3)$;

(8)\quad $\L(1,0) \otimes \L(2,2) \simeq \T(3,2)$;

(9)\quad $\L(2,0) \otimes \L(2,0) \simeq \T(4,0) \oplus \T(2,1)$;

(10)\quad $\L(2,0) \otimes \L(1,1) \simeq \T(3,1) \oplus \T(0,1)$;

(11)\quad $\L(2,0) \otimes \L(0,2) \simeq \T(2,2) \oplus \T(1,1)$;

(12)\quad $\L(2,0) \otimes \L(2,1) \simeq \T(4,1) \oplus \T(2,2)$;

(13)\quad $\L(2,0) \otimes \L(1,2) \simeq \T(3,2)\oplus \T(0,2) \oplus
  \T(1,0)$;

(14)\quad $\L(2,0) \otimes \L(2,2) \simeq \T(4,2) \oplus \T(2,3)$;

(15)\quad $\L(1,1) \otimes \L(1,1) \simeq \T(2,2) \oplus \T(0,0) \oplus
\M$; 

(16)\quad $\L(1,1) \otimes \L(2,1) \simeq \T(3,2) \oplus \T(4,0)
\oplus \T(1,0)$;

(17)\quad $\L(1,1) \otimes \L(2,2) \simeq \T(3,3) \oplus \T(2,2)$;

(18)\quad $\L(2,1) \otimes \L(2,1) \simeq \T(4,2) \oplus \T(5,0)
\oplus \T(2,3) \oplus \T(3,1)$;

(19)\quad $\L(2,1) \otimes \L(1,2) \simeq \T(3,3) \oplus 2\T(2,2) \oplus
 \T(1,1)$;

(20)\quad $\L(2,1) \otimes \L(2,2) \simeq \T(4,3) \oplus 2\T(3,2)
 \oplus \T(2,4)$;

(21)\quad $\L(2,2) \otimes \L(2,2) \simeq \T(4,4) \oplus \T(3,3)
 \oplus \T(5,2) \oplus \T(2,5) \oplus 3 \T(2,2)$.}

\noindent
Thus the family $\F$ is in this case given by the twenty-five tilting
modules $\{ \T(a,b): 0 \le a,b \le 4 \}$ along with the six
``exceptional'' modules $$\{ \T(5,0), \T(0,5), \T(5,2), \T(2,5),
\L(1,1), \M \}.$$ All members of $\F$ except $\L(1,1)$ and $\M$ are
tilting modules.

\item The uniserial members of $\F$ have the following structure: {\rm

$\T(0,0) = [(0,0)]$;\quad $\T(1,0) = [(1,0)]$;\quad $\T(2,0) =
  [(2,0)]$;

$\T(1,1) = [(0,0), (1,1), (0,0)]$; \quad $\T(2,1) = [(2,1)]$;

$\T(4,0) = [(0,2), (4,0), (0,2)]$; \quad $\T(3,1) = [(1,2), (3,1),
    (1,2)]$;

$\T(2,2) = [(2,2)]$; \quad $\T(5,0) = [(0,1), (5,0), (0,1)]$;\quad
  $\T(5,2) = [(5,2)]$.  }

\noindent
The structure of the non-uniserial rigid members of $\F$ is given
below (symmetric cases omitted):
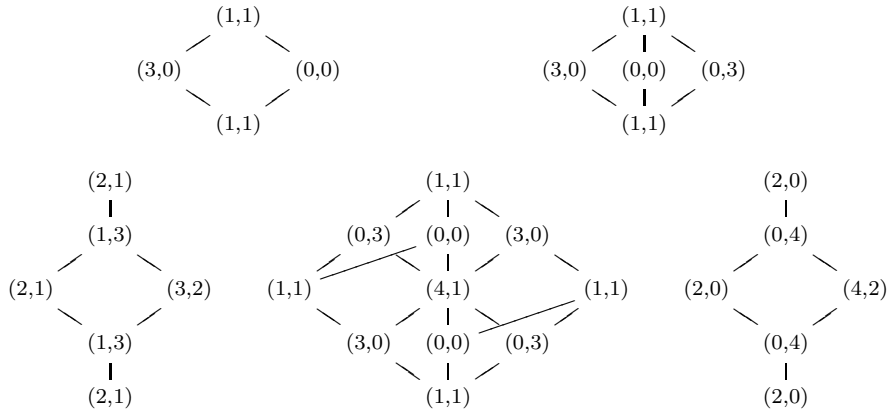
\begin{center}
\begin{tabular*}{0.57\textwidth}{@{\extracolsep{\fill}}lr}
\def\objectstyle{\scriptstyle}
\xymatrix@=6pt{
  & (1,1) \ar@{-}[dl] \ar@{-}[dr] & \\
  (3,0) & & (0,0) \\
  & (1,1) \ar@{-}[ul] \ar@{-}[ur] &
} 
&
\def\objectstyle{\scriptstyle}
\xymatrix@=6pt{
 & (1,1) \ar@{-}[dr] \ar@{-}[dl] & \\
  (3,0) & (0,0)\ar@{-}[u] \ar@{-}[d]   & (0,3) \\
 & (1,1) \ar@{-}[ur] \ar@{-}[ul] &
}
\end{tabular*}
\end{center}

\begin{center}
\begin{tabular}{ccc}
\def\objectstyle{\scriptstyle}
\xymatrix@=6pt{
  & (2,1) \ar@{-}[d] & \\
  & (1,3) \ar@{-}[dl] \ar@{-}[dr] & \\
(2,1) \ar@{-}[dr] & & (3,2) \ar@{-}[dl] \\
  & (1,3) \ar@{-}[d] & \\
  & (2,1) & 
}
&
\def\objectstyle{\scriptstyle}
\xymatrix@=6pt{
  & & (1,1) \ar@{-}[d] \ar@{-}[dl] \ar@{-}[dr]& & \\
  & (0,3)\ar@{-}[dl]\ar@{-}[dr] & (0,0)\ar@{-}[d]\ar@{-}[dll] 
   & (3,0)\ar@{-}[dl] \ar@{-}[dr] & \\
(1,1) \ar@{-}[dr] & & (4,1) \ar@{-}[dl] \ar@{-}[d]\ar@{-}[dr] 
  & & (1,1)\ar@{-}[dl]\ar@{-}[dll] \\
  & (3,0) & (0,0)  & (0,3) & \\
   & & (1,1) \ar@{-}[u] \ar@{-}[ul] \ar@{-}[ur]& & 
}
&
\def\objectstyle{\scriptstyle}
\xymatrix@=6pt{
  & (2,0) \ar@{-}[d] & \\
  & (0,4) \ar@{-}[dl] \ar@{-}[dr] & \\
(2,0) \ar@{-}[dr] & & (4,2) \ar@{-}[dl] \\
  & (0,4) \ar@{-}[d] & \\
  & (2,0) & 
}
\end{tabular}
\end{center}
all of which are tilting modules excepting the module $\M$ (which does
not have a highest weight) pictured at the upper right. Finally, there
are four members of $\F$, namely $\T(3,3)$, $\T(4,3)$, $\T(3,4)$, and
$\T(4,4)$, whose structure is not rigid, which are not
pictured. Analysis of their structure requires other methods (see the
Appendix).

\item Each member of $\F$ except $\T(5,2) = \L(5,2)$, $\T(2,5) =
  \L(2,5)$ has simple $G_1T$-socle (and head) of highest
  weight belonging to $X_1$.
\end{enumerate}
\end{prop*}

\begin{rmk}
The Alperin diagram for $\T(4,1)$ given above is one of several
possibilities. When a module has a direct sum of two or more copies of
the same simple on a given socle layer, there may be more than one
diagram.
\end{rmk}

The proof of the proposition will be given in
\ref{ss:pfp=3}--\ref{partc}. First we consider some consequences and
look at a few examples.

\begin{cor*}
  Let $p=3$. Given arbitrary dominant weights $\lambda, \mu$ write
  $\lambda = \sum \lambda^j p^j$, $\mu = \sum \mu^j p^j$ with each
  $\lambda^j, \mu^j \in X_1$.
  \begin{enumerate}[label={\rm(\alph*)},leftmargin=*,itemsep=0.5em]
  \item In the decomposition \eqref{eq:star}, each term in the direct
    sum not involving a tensor factor of the form $\T(5,2)$, $\T(2,5)$
    is indecomposable.

  \item $\L(\lambda) \otimes \L(\mu)$ is indecomposable if and only if
    for each $j \ge 0$ the unordered pair $\{\lambda^j, \mu^j\}$ is
    one of the cases $\{(1,0), (0,1)]$, $\{(1,0), (2,0)]$, $\{(1,0),
        (2,2)]$, $\{(0,1), (0,2)]$, $\{(0,1), (2,2)]$ or one of
              $\lambda^j$, $\mu^j$ is the zero weight $(0,0)$.

  \item Let $m$ be the maximum $j$ such that at least one of
    $\lambda^j$, $\mu^j$ is non-zero. Then $\L(\lambda) \otimes
    \L(\mu)$ is indecomposable tilting, isomorphic to
    $\T(\lambda+\mu)$, if and only if: {\rm(i)} for each $0\le j \le
    m-1$, one of $\lambda^j, \mu^j$ is minuscule and the other is the
    Steinberg weight $(2,2)$, and {\rm(ii)} $\{\lambda^m, \mu^m\}$ is
    one of the cases listed in part {\rm(b)}.
  \end{enumerate}
\end{cor*}

\begin{proof}
  The proof is entirely similar to the proof of the corresponding
  result in the $p=2$ case. We leave the details to the reader.
\end{proof}

\begin{examples}
  (i) We work out the indecomposable direct summands of $\L(5,4)
  \otimes \L(4,5)$, using information from part (a) of the proposition
  and following the procedure of Section \ref{intro:gen}:
  \begin{align*}
    \L(5,4) &\otimes \L(4,5) \\
    & \simeq \big(\L(2,1) \otimes
    \L(1,1)^{[1]}\big) \otimes \big(\L(1,2) \otimes \L(1,1)^{[1]}\big)
    \\ &\simeq \big(\L(2,1)\otimes \L(1,2)\big) \otimes \big(\L(1,1)
    \otimes \L(1,1) \big)^{[1]} \\
    &\simeq \big( \T(3,3) \oplus 2\T(2,2) \oplus \T(1,1) \big) \otimes 
    \big( \T(2,2) \oplus \T(0,0) \oplus M \big)^{[1]} \\
    & \simeq \big(\T(3,3) \otimes \T(2,2)^{[1]}\big) \oplus 
    \big(\T(3,3) \otimes \T(0,0)^{[1]}\big) \oplus
    \big(\T(3,3) \otimes \M^{[1]}\big) \\
    & \qquad \oplus 2 \big(\T(2,2) \otimes \T(2,2)^{[1]}\big) \oplus
    2 \big(\T(2,2) \otimes \T(0,0)^{[1]}\big) \oplus
    2 \big(\T(2,2) \otimes \M^{[1]}\big) \\
    & \qquad\oplus\big(\T(1,1) \otimes \T(2,2)^{[1]}\big) \oplus
    \big(\T(1,1) \otimes \T(0,0)^{[1]}\big) \oplus
    \big(\T(1,1) \otimes \M^{[1]}\big)\\
    &\simeq \T(9,9) \oplus \T(3,3) \oplus 
     \big(\T(3,3) \otimes \M^{[1]}\big) \oplus 2\T(8,8) \oplus 2\T(2,2)\\
     & \qquad\oplus 2 \big(\T(2,2) \otimes \M^{[1]}\big) \oplus \T(7,7) 
     \oplus \T(1,1) \oplus \big(\T(1,1) \otimes \M^{[1]}\big).
  \end{align*}
  We applied \eqref{eq:tptilt} to get the last line of the
  calculation. 

  \noindent
  (ii) Illustrating part (b) of the corollary we have $\L(3,1) \otimes
  \L(1,3) \simeq \L(0,1) \otimes \L(1,0)^{[1]} \otimes \L(1,0) \otimes
  \L(0,1)^{[1]} \simeq \T(1,1) \otimes \T(1,1)^{[1]}$ , which is
  indecomposable but not tilting. 

  \noindent
  (iii) To illustrate part (c) of the corollary we have for instance
  $\L(4,0) \otimes \L(8,8) \simeq \T(12,8)$ or $\L(5,2) \otimes
  \L(5,4) \simeq \T(10,6)$.
\end{examples}

\subsection{}\label{p3:decisive}
We now discuss the problem of computing the indecomposable direct
summands (and their multiplicities) of $\L(\lambda) \otimes \L(\mu)$
for arbitrary $\lambda, \mu \in X^+$, in the more difficult case where
a direct summand on the right hand side of \eqref{eq:star} is not
necessarily indecomposable. 

It will be convenient to introduce the notation $\F_0$ for the set $\F
- \{ \T(5,2), \T(2,5) \}$. Then Corollary \ref{main:3}(a) says that a
direct summand $S = \bigotimes_{j \ge 0}\, I_j^{[j]}$ in
\eqref{eq:star} is indecomposable whenever all its tensor factors
$I_j$ belong to $\F_0$.

Consider a summand $S = \bigotimes_{j \ge 0}\, I_j^{[j]}$ in
\eqref{eq:star} which is possibly not indecomposable. By Corollary
\ref{main:3}(a), such a summand must have one or more tensor
multiplicands of the form $\T(5,2)$ or $\T(2,5)$. Suppose that in the
summand in question $I_k$ is $\T(5,2)$ or $\T(2,5)$. We use the fact
that $\T(5,2) = \L(5,2) \simeq \L(2,2) \otimes \L(1,0)^{[1]}$, and
similarly $\T(2,5) = \L(2,5) \simeq \L(2,2) \otimes
\L(0,1)^{[1]}$. Thus we are forced to consider the possible splitting
of $\L(1,0) \otimes I_{k+1}$ or $\L(0,1) \otimes I_{k+1}$ in `degree'
$k+1$. (By `degree' here we just mean the level of $j$ in the twisted
tensor product occurring in a direct summand of the right-hand-side of
\eqref{eq:star}.)  There are two cases.

We consider first the case where $I_{k+1}$ is not tilting, i.e.,
$I_{k+1}$ is either $\L(1,1)$ or $\M$. So we need to split $\L(1,0)
\otimes \L(1,1)$, $\L(0,1) \otimes \L(1,1)$, $\L(1,0) \otimes \M$ or
$\L(0,1) \otimes \M$.  The first two cases are already covered by
Corollary \ref{main:3}(a), so we just need to consider the last two.
But a simple calculation with characters and consideration of linkage
classes shows that
\begin{equation}
  \begin{gathered}
  \L(1,0) \otimes \M \simeq \T(3,1) \oplus \T(1,0) \oplus \T(4,0);\\ 
  \L(0,1) \otimes \M \simeq \T(1,3) \oplus \T(0,1) \oplus \T(0,4)
  \end{gathered}
\end{equation}
and the summands are once again members of $\F$ with restricted
socles, so these cases present no problem.

We are left with the case where $I_{k+1}$ is tilting. Then this
splitting can be computed since $\L(1,0) = \T (1,0)$ and $\L(0,1) =
\T(0,1)$ are tilting, so we are just splitting a tensor product of two
tilting modules into a direct sum of indecomposable tilting modules,
which can always be done. This new decomposition produces only tilting
modules in the family $\F$ except when $I_{k+1}$ is one of the
following cases:
\[
\T(5,0), \T(4,1), \T(4,2), \T(5,2), \T(4,3), \text{ and } \T(4,4)
\]
or one of their symmetric cousins. Up to lower order terms which again
belong to $\F$, these possibilities, when tensored by $\L(1,0)$ or
$\L(0,1)$, produce the new tilting modules 
\begin{equation}
\T(6,0), \T(5,1), \T(6,2), \T(5,3), \text{ and } \T(5,4)
\end{equation}
and of course their symmetric versions. Now by Donkin's tensor product
theorem we have a twisted tensor product decomposition for the last
three of these, in terms of members of $\F$: 
\begin{equation}
\begin{aligned}
  \T(6,2) &\simeq \T(3,2) \otimes \T(1,0)^{[1]},\\ 
  \T(5,3) &\simeq \T(2,3) \otimes \T(1,0)^{[1]},\\
  \T(5,4) &\simeq \T(2,4) \otimes \T(1,0)^{[1]}.  
\end{aligned}
\end{equation}
Hence, those summands and their symmetric versions present no
problem. Finally, if $\T(5,0)$ or $\T(4,1)$ is tensored by $\L(1,0)$
then, modulo lower order terms which belong to $\F$, we obtain the new
summands $\T(6,0)$ and $\T(5,1)$ which are not members of $\F$ and do
not admit a twisted tensor product decomposition. However, these
summands must have simple restricted $G_1T$-socles, since they are
embedded in $\T(4,4)$ and $\T(4,3)$, respectively. This is shown by
translation arguments, similar to those in \ref{partc} ahead. Thus we
have proved the following result.

\begin{prop*} Let $p=3$ and $G=\SL_3$.
\begin{enumerate}[label={\rm(\alph*)},leftmargin=*,itemsep=0.5em]
\item  Any tensor product of the form $\L(1,0) \otimes I$ or $\L(0,1)
  \otimes I$, where $I$ is an indecomposable tilting module in $\F$,
  is expressible as a twisted tensor product of modules which are
  either tilting modules in $\F_0$ or are one of the ``extra'' modules
  $\T(6,0)$, $\T(5,1)$, $\T(0,6)$ or $\T(1,5)$.

\item The extra modules have simple restricted $G_1 T$-socles.

\item For general $\lambda, \mu \in X^+$, the indecomposable direct
  summands of $\L(\lambda) \otimes \L(\mu)$ are all expressible as
  twisted tensor products of modules from the family
\begin{equation*}
  \begin{aligned}
  \F' &= \F_0 \cup \{\T(6,0), \T(5,1), \T(0,6), \T(1,5) \}\\ &=
  \{\T(a,b)\colon 0\le a,b \le 4\}\ \cup \\ &\qquad \{\T(6,0), \T(5,1),
  \T(5,0), \T(0,6), \T(1,5), \T(0,5), \L(1,1), \M \}.
  \end{aligned}
\end{equation*}
Note that all members of $\F'$ have simple $G_1T$-socle of restricted
highest weight.
\end{enumerate}
\end{prop*}

\begin{example}
  \noindent
  We consider an example where the direct summands on the
  right hand side of \eqref{eq:star} are not all
  indecomposable:
  \begin{align*}
    \L(2,2) &\otimes \L(5,2) \simeq \L(2,2) \otimes \L(2,2) \otimes
    \L(1,0)^{[1]}\\ & \simeq \big( \T(4,4) \oplus \T(3,3) \oplus
    \T(5,2) \oplus \T(2,5) \oplus 3\T(2,2) \big) \otimes \L(1,0)^{[1]}\\
    & \simeq \T(7,4) \oplus \T(6,3) \oplus \T(8,2) \oplus \T(2,5) \oplus
    \T(5,5) \oplus 3 \T(5,2).
  \end{align*}
  The second line comes from equation (21) in Proposition
  \ref{main:3}(a), and to get the last line one applies
  \eqref{eq:tptilt} repeatedly, using Proposition \ref{main:3}(a)
  again as needed. For instance, using equation (1) from Proposition
  \ref{main:3}(a) we have
  \begin{align*}
   \T(5,2) \otimes \L(1,0)^{[1]} &\simeq \L(2,2) \otimes
   \big(\L(1,0) \otimes \L(1,0)\big)^{[1]}\\ & \simeq
   \L(2,2) \otimes \big(\T(2,0) \oplus \T(0,1) \big)^{[1]} \\
   & \simeq \T(8,2) \oplus \T(2,5)
  \end{align*}
  and using equation (2) from Proposition \ref{main:3}(a) we have
  \begin{align*}
    \T(2,5) \otimes \L(1,0)^{[1]} &\simeq \L(2,2) \otimes \big(\L(0,1)
    \otimes \L(1,0)\big)^{[1]} \\
    & \simeq \L(2,2) \otimes \T(1,1)^{[1]} \simeq \T(5,5).
  \end{align*} 
\end{example}

\subsection{}\label{ss:pfp=3}
We now embark upon the proof of Proposition \ref{main:3}.  First we
compute the composition factor multiplicities of the restricted tensor
products. (Recall that $\chi_p(\lambda) = \L(\lambda)$ is the formal
character of $\L(\lambda)$.)
\begin{enumerate}
\makeatletter
\renewcommand{\labelenumi}{(\theenumi) }
\makeatother
\item $\chi_p(1,0)\cdot \chi_p(1,0) = \chi_p(2,0) + \chi_p(0,1)$;

\item $\chi_p(1,0)\cdot \chi_p(0,1) = \chi_p(1,1) + 2\chi_p(0,0)$;

\item $\chi_p(1,0)\cdot \chi_p(2,0) = \chi_p(3,0) + 2\chi_p(1,1) + \chi_p(0,0)$;

\item $\chi_p(1,0)\cdot \chi_p(1,1) = \chi_p(2,1) + \chi_p(0,2)$;

\item $\chi_p(1,0)\cdot \chi_p(0,2) = \chi_p(1,2) + \chi_p(0,1)$;

\item $\chi_p(1,0)\cdot \chi_p(2,1) = \chi_p(3,1) + 2\chi_p(1,2) + \chi_p(2,0)$;

\item $\chi_p(1,0)\cdot \chi_p(1,2) = \chi_p(2,2) + \chi_p(0,3) +
2\chi_p(1,1) + \chi_p(0,0)$;

\item $\chi_p(1,0)\cdot \chi_p(2,2) = \chi_p(3,2) + 2\chi_p(1,3) + 3\chi_p(2,1)$;

\item $\chi_p(2,0)\cdot \chi_p(2,0) = \chi_p(4,0) + \chi_p(2,1) + 2\chi_p(0,2)$;

\item $\chi_p(2,0)\cdot \chi_p(1,1) = \chi_p(3,1) + 2\chi_p(1,2) + \chi_p(0,1)$;

\item $\chi_p(2,0)\cdot \chi_p(0,2) = \chi_p(2,2) + \chi_p(1,1) + 2\chi_p(0,0)$;

\item $\chi_p(2,0)\cdot \chi_p(2,1) = \chi_p(4,1) + \chi_p(2,2) +
2\chi_p(0,3) + 2\chi_p(3,0) + 4\chi_p(1,1) + 2\chi_p(0,0)$;

\item $\chi_p(2,0)\cdot \chi_p(1,2) = \chi_p(3,2) + 2\chi_p(1,3) +
3\chi_p(2,1) + \chi_p(0,2) + \chi_p(1,0)$;

\item $\chi_p(2,0)\cdot \chi_p(2,2) = \chi_p(4,2) + \chi_p(2,3) +
2\chi_p(0,4) + 2\chi_p(3,1) + 3\chi_p(1,2) + 3\chi_p(2,0)$;

\item $\chi_p(1,1)\cdot \chi_p(1,1) = \chi_p(2,2) + \chi_p(0,3) +
\chi_p(3,0) + 2\chi_p(1,1) + 2\chi_p(0,0)$;

\item $\chi_p(1,1)\cdot \chi_p(2,1) = \chi_p(3,2) + 2\chi_p(1,3) +
\chi_p(4,0) + 3\chi_p(2,1) + 2\chi_p(0,2) + \chi_p(1,0)$;

\item $\chi_p(1,1)\cdot \chi_p(2,2) = \chi_p(3,3) + 2\chi_p(1,4) +
2\chi_p(4,1) + \chi_p(2,2) + 4\chi_p(0,3) + 4\chi_p(3,0) +
6\chi_p(1,1) + 5\chi_p(0,0)$;

\item $\chi_p(2,1)\cdot \chi_p(2,1) = \chi_p(4,2) + \chi_p(2,3) +
2\chi_p(0,4) + \chi_p(5,0) + 3\chi_p(3,1) + 5\chi_p(1,2) +
3\chi_p(2,0) + 2\chi_p(0,1)$;

\item $\chi_p(2,1)\cdot \chi_p(1,2) = \chi_p(3,3) + 2\chi_p(1,4) +
2\chi_p(4,1) + 2\chi_p(2,2) + 4\chi_p(0,3) + 4\chi_p(3,0) +
7\chi_p(1,1) + 7\chi_p(0,0)$;

\item $\chi_p(2,1)\cdot \chi_p(2,2) = \chi_p(4,3) + \chi_p(2,4) +
2\chi_p(0,5) + 2\chi_p(5,1) + 2\chi_p(3,2) + 4\chi_p(1,3) +
2\chi_p(4,0) + 6\chi_p(2,1) + 3\chi_p(0,2) + 5\chi_p(1,0)$;

\item $\chi_p(2,2)\cdot \chi_p(2,2) = \chi_p(4,4) + \chi_p(2,5) +
2\chi_p(0,6) + \chi_p(5,2) + 3\chi_p(3,3) + 6\chi_p(1,4) +
2\chi_p(6,0) + 6\chi_p(4,1) + 3\chi_p(2,2) + 8\chi_p(0,3) +
8\chi_p(3,0) + 11\chi_p(1,1) + 15\chi_p(0,0)$.
\end{enumerate}

It is important to proceed inductively through the cases, so that the
structure of smaller tilting modules is available by the time the
argument reaches the higher, more complicated, cases. We order the
cases as listed in part (a) of Proposition \ref{main:3}. In each case,
one starts by partitioning the composition factors into blocks.  This
amounts to looking at linkage classes of the highest weights of those
composition factors.

In cases (1)--(9), (11)--(14) the argument is entirely similar to the
arguments already used in characteristic 2, in the proof of
Proposition \ref{ss:rtpd}. In these cases we know that the tensor
product in question is tilting, and it turns out that each linkage
class determines an indecomposable direct summand. This uses the
contravariant self-duality of the tilting modules and the structural
information in \ref{Weyl:p3} for the Weyl modules, which forces a
lower bound on the composition length of the tilting module in
question, and it turns out that this lower bound agrees with the upper
bound provided by the linkage class. 

As an example, let us examine the argument in the case (12), for the
tensor product $\L(2,0) \otimes \L(2,1)$. The linkage classes are
\[
\{(2,2)\} \cup \{(4,1), (3,0), (0,3), (1,1), (0,0) \}. 
\]
By highest weight considerations, we must have a single copy of
$\T(4,1)$ in $\L(2,0) \otimes \L(2,1)$. Linkage forces a copy of
$\L(2,2) = \T(2,2)$ to split off as well. Now $\T(4,1)$ has a
submodule isomorphic to $\Delta(4,1)$, so $\L(1,1)$ is contained in
its socle. By self-duality of $\T(4,1)$, we must have another copy of
$\L(1,1)$ in the top of $\T(4,1)$, so we are forced to put a copy of
$\Delta(1,1)$ at the top of $\T(4,1)$. Looking at the structure of
$\Delta(4,1)$ in 5.1.6, we see that $\T(4,1)$ must also have at least
one copy of $\Delta(3,0)$ and $\Delta(0,3)$ in its
$\Delta$-filtration. At this point we are finished, since this
accounts for all available composition factors (with their
multiplicities) from the linkage class, so we conclude that $\L(2,0)
\otimes \L(2,1) \simeq \T(4,1) \oplus \T(2,2)$. The structure of
$\T(4,1)$ is nearly forced, because of its self-duality, the fact that
all the Ext groups between simple factors is known, and the fact that
$\T(4,1)$ must have both $\Delta$ and $\nabla$-filtrations. In
\ref{partc} we will see that $\T(0,3)$ is isomorphic to a submodule
of $\T(4,1)$, which finishes the determination of the structure of
$\T(4,1)$.

\subsection{}\label{lem:L11}
In case (10), one cannot immediately conclude that $\L(2,0) \otimes
\L(1,1)$ is tilting since $\L(1,1)$ is not tilting, so we must proceed
differently. However, we observe the following, which immediately
implies that in fact our tensor product is tilting.

\begin{lem*}
Let $V$ be a simple Weyl module and let $\Delta(\lambda)$ be a Weyl
module of highest weight $\lambda$.  If the composition factors of $V
\otimes \operatorname{rad} \Delta(\lambda)$ and $V \otimes L(\lambda)$
lie in disjoint blocks, then $V \otimes L(\lambda)$ is tilting.
\end{lem*}

\begin{proof}
$V \otimes \Delta(\lambda)$ has a $\Delta$-filtration, by the
  Wang--Donkin--Mathieu result (see \cite[II.4.21]{Jantzen}). Now as
  $V \otimes \operatorname{rad} \Delta(\lambda)$ and $V \otimes
  L(\lambda)$ have no common linkage classes there can be no
  non-trivial extensions between these modules, by the linkage
  principle.  Thus $V \otimes \Delta(\lambda) = \big(V \otimes
  \operatorname{rad} \Delta(\lambda)\big) \oplus \big(V \otimes
  L(\lambda)\big)$. As $V \otimes \Delta(\lambda)$ has a
  $\Delta$-filtration this implies $V \otimes L(\lambda)$ does
  also. As it is the tensor product of two simple (therefore
  contravariantly self dual) modules it is itself contravariantly self
  dual and so has a $\nabla$-filtration. Therefore it is tilting.
\end{proof}

Now we may proceed as usual. Looking at the character of $\L(2,0)
\otimes \L(1,1)$ we find that there are two linkage classes for the
highest weights of the composition factors, namely $\{(0,1)\}$ and
$\{(3,1), (1,2)\}$. Since the multiplicity of $\L(0,1)$ is 1, it must
give a simple tilting summand $\T(0,1)$. Now $\T(3,1)$ must be a
summand by highest weight consideration, and the usual argument forces
it to have at least composition length three, which forces equality of
the upper and lower bounds, so the structure is $\T(3,1) = [(1,2),
  (3,1), (1,2)]$ and we have $\L(2,0) \otimes \L(1,1) \simeq \T(3,1)
\oplus \T(0,1)$. This takes care of case (10) in our list. 

Case (16) follows similarly, making use again of the above lemma to
conclude that $\L(1,1) \otimes \L(2,1)$ is tilting. We note that at
this stage we may assume that the structure of $\T(3,2)$ and $\T(4,0)$
are already known, since they come up in the earlier cases (8),
(9). So one easily concludes from this and the linkage classes that
$\L(1,1) \otimes \L(2,1) \simeq \T(3,2) \oplus \T(4,0) \oplus
\T(1,0)$.

\subsection{} \label{M}
We now consider case (15). Since $\L(1,1)$ is not tilting, it is
unclear whether or not $\L(1,1) \otimes \L(1,1)$ is tilting. In fact
it is not, and analysis of this case is more difficult. First, looking
at the character and the linkage classes (there are two) we observe
that a copy of the Steinberg module $\T(2,2)=\L(2,2)$ splits off as a
direct summand. The remaining composition factors of the tensor
product all lie in the same linkage class, but it turns out that a
copy of the trivial module splits off, as we show below.

From properties of duals and previous calculations it follows that
\[ \tag{1}
\begin{aligned}
\dim_K \Hom_G&(\L(0,0), \L(1,1)\otimes \L(1,1))\\ &= \dim_K
\Hom_G(\L(0,0) \otimes \L(1,1), \L(1,1))\\ &= \dim_K \Hom_G(\L(1,1),
\L(1,1)) = 1;
\end{aligned}
\]

\[ \tag{2}
\begin{aligned}
\dim_K \Hom_G&(\T(1,1), \L(1,1)\otimes \L(1,1))\\ &= \dim_K
\Hom_G(\L(1,0) \otimes \L(0,1), \L(1,1) \otimes \L(1,1)) \\ &= \dim_K
\Hom_G(\L(1,1) \otimes \L(0,1), \L(1,1) \otimes \L(0,1))\\ &= \dim_K
\Hom_G(\L(1,2) \oplus \L(2,0), \L(1,2) \oplus \L(2,0)) = 2;
\end{aligned}
\]

\[ \tag{3}
\begin{aligned}
\dim_K \Hom_G&(\L(3,0), \L(1,1)\otimes \L(1,1))\\ &= \dim_K
\Hom_G(\L(1,1) \otimes \L(3,0), \L(1,1))\\ &= \dim_K \Hom_G(\L(4,1),
\L(1,1)) = 0
\end{aligned}
\]
and, by symmetry, an equality similar to (3) holds, in which $(3,0)$ is
replaced by $(0,3)$. We also observe that
\[ \tag{4}
  \Hom_G(\L(1,1), \L(1,1)\otimes \L(1,1)) \simeq \Hom_G(\L(1,1)\otimes
  \L(1,1), \L(1,1))
\]
By (1), (3), and (4) we see that the socle of $\L(1,1) \otimes
\L(1,1)$ is either: (a) $\L(2,2) \oplus \L(0,0)$, or (b) $\L(2,2)
\oplus \L(0,0) \oplus \L(1,1)$. 

From the structure of the Weyl modules in question we know (see
e.g.\ \cite[II.4.14]{Jantzen}) all the $\Ext^1$ groups between the
simple modules of interest here. Combining this with self-duality
would force the structure of the non-simple direct summand of $\L(1,1)
\otimes \L(1,1)$ to be given by one of the following diagrams:
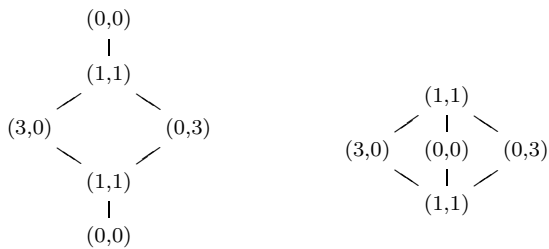
\[
\def\objectstyle{\scriptstyle}
\xymatrix@=6pt{
 & (0,0) \ar@{-}[d] \\
 & (1,1) \ar@{-}[dr] \ar@{-}[dl] & \\
  (3,0) &   & (0,3) \\
 & (1,1) \ar@{-}[ur] \ar@{-}[ul] & \\
 & (0,0) \ar@{-}[u]
}
\qquad \qquad 
\def\objectstyle{\scriptstyle}
\xymatrix@=6pt{ {\ } \\ \\
 & (1,1) \ar@{-}[dr] \ar@{-}[dl] & \\
  (3,0) & (0,0)\ar@{-}[u] \ar@{-}[d]   & (0,3) \\
 & (1,1) \ar@{-}[ur] \ar@{-}[ul] &
}
\]
where the left diagram corresponds with possibility (a) and the right
with possibility (b).  However, the left diagram would contradict
(2). Hence, possibility (a) is in fact ruled out, and we are left with
possibility (b). It follows that $\L(1,1) \otimes \L(1,1) \simeq
\L(2,2) \oplus \L(0,0) \oplus \M$, as claimed.

\subsection{}\label{T33}
There are just five cases remaining in the proof of Proposition
\ref{main:3}, namely cases (17)--(21).  We now consider case (17). The
module $\L(1,1) \otimes \L(2,2)$ is tilting by Lemma \ref{lem:L11}, so
by highest weight considerations $\T(3,3)$ is a direct summand. This
is also justified by Pillen's Theorem (see \ref{ss:genPillen}). The
character of $\T(3,3)$ may be computed by \eqref{eq:Donkinch2}, which
shows that it has a $\Delta$-filtration with $\Delta$-factors
isomorphic to
\[
 \Delta(3,3),\ \Delta(4,1),\ \Delta(1,4),\ \Delta(3,0),\
 \Delta(0,3),\ \Delta(1,1)
\]
each occurring with multiplicity one. This accounts for all the
composition factors appearing in the character of $\L(1,1) \otimes
\L(2,2)$, except for one copy of the Steinberg module $\T(2,2) =
\L(2,2)$. Hence we conclude that
\[
 \L(1,1) \otimes \L(2,2) \simeq \T(3,3) \oplus \T(2,2). 
\]

\subsection{}
$\L(2,1) \otimes \L(2,1)$ is tilting since $\L(2,1)$ is, so by highest
weight considerations a copy of $\T(4,2)$ splits off as a direct
summand. The structure of $\T(4,2)$ was determined in a previous case
of the proof. Subtracting its character from the character of $\L(2,1)
\otimes \L(2,1)$, we see that the highest weight of what remains is
$(5,0)$, so a copy of $\T(5,0)$ must split off as well. The linkage
class of $(5,0)$ contains only two weights $\{(5,0), (0,1)\}$ and from
this and the known structure of the Weyl modules it follows easily
that $\T(5,0)$ is uniserial with structure $\T(5,0) = [(0,1), (5,0),
  (0,1)]$.  Now highest weight and character considerations force the
remaining summands to be one copy of $\T(2,3)$ and one copy of
$\T(3,1)$. Hence
\[
\L(2,1) \otimes \L(2,1) \simeq \T(4,2) \oplus \T(5,0) \oplus \T(2,3)
\oplus \T(3,1).
\]
We note we can assume that $\T(3,1)$ and $\T(2,3)$ are known at this
point, since they arise in earlier cases of the proof. (Actually, to
be precise $\T(2,3)$ doesn't arise in any earlier case, but its
symmetric cousin $\T(3,2)$ does.)

\subsection{}
$\L(2,1) \otimes \L(1,2)$ is tilting since both $\L(2,1)$ and
$\L(1,2)$ are, so by highest weight considerations a copy of $\T(3,3)$
splits off as a direct summand. The character of $\T(3,3)$ was
computed already in \ref{T33}, so by character considerations one
easily deduces that
\[
 \L(2,1) \otimes \L(1,2) \simeq \T(3,3) \oplus 2\T(2,2) \oplus
 \T(1,1).
\]
Of course, the character of $\T(1,1)$ is already known by an earlier
case of the proof.

\subsection{}\label{T43}
$\L(2,1) \otimes \L(2,2)$ is tilting since both $\L(2,1)$ and
$\L(2,2)$ are, so by highest weight considerations a copy of $\T(4,3)$
splits off as a direct summand. From \cite[\S2.1, Lemma
  5]{Donkin:coho} we compute its $\Delta$-factors to be
\[
  \Delta(4,3),\ \Delta(5,1),\ \Delta(0,5),\ \Delta(1,0).
\]
One sees also that $\T(4,3)$ has simple socle of highest weight
$(1,0)$ by arguments similar to those in \ref{T33}. From character
computations one now shows that
\[
 \L(2,1) \otimes \L(2,2) \simeq \T(4,3) \oplus 2\T(3,2) \oplus \T(2,4).
\]
The structure of $\T(3,2)$ is available by a previous case of the
proof, and the structure of $\T(2,4)$ follows by symmetry from that of
$\T(4,2)$, again a previous case.

\subsection{}\label{T44}
$\L(2,2) \otimes \L(2,2)$ is tilting since $\L(2,2)$ is, so by highest
weight considerations a copy of $\T(4,4)$ must split off as a direct
summand. The $\Delta$-factor multiplicities of $\T(4,4)$ are computed
by \cite[\S2.1, Lemma 5]{Donkin:coho} to
be
\[ 
 \Delta(4,4),\ \Delta(6,0),\ \Delta(0,6),\ \Delta(3,3),\
 \Delta(4,1),\ \Delta(1,4),\ \Delta(1,1),\ \Delta(0,0)
\]
each of multiplicity one. From this, using the character of $\L(2,2)
\otimes \L(2,2)$ it follows by highest weight considerations, after
subtracting the character of $\T(4,4)$, that a copy of $\T(3,3)$ must
also split off as a direct summand. Then it easily follows that
\[
 \L(2,2) \otimes \L(2,2) \simeq \T(4,4) \oplus \T(3,3) \oplus \T(5,2)
 \oplus \T(2,5) \oplus 3 \T(2,2)
\]
where $\T(5,2) = \L(5,2)$, $\T(2,5)= \L(2,5)$, and $\T(2,2) =
\L(2,2)$. 

At this point the proof of Proposition \ref{main:3}(a), (b) is
complete. 

\subsection{}\label{partc}
It remains to prove the claim in part (c) of Proposition \ref{main:3}.
It is known that Donkin's conjecture holds for $\SL_3$, as discussed
at the beginning of \ref{gen:tp-thm}, so $\T((p-1)\rho+\lambda)$ is
as a $G_1T$-module isomorphic to $\hatQ_1((p-1)\rho + w_0 \lambda)$
for any $\lambda \in X_1$.  Thus $\T(a,b)$ has simple $G_1T$-socle of
restricted highest weight, for any $2 \le a, b \le 4$.  Moreover, the
claim is true of $\T(0,0)$, $\T(1,0)$, $\T(2,0)$, $\T(2,1)$, $\L(1,1)$
and their symmetric counterparts, since these are all simple
$G$-modules of restricted highest weight.

For $\lambda = (1,1)$ and $(5,0)$ one easily checks by direct
computation that $\Delta(\lambda)$, which is a non-split extension
between two simple $G$-modules, remains non-split upon restriction to
$G_1T$. It then follows that $\T(\lambda)$ has simple $G_1T$-socle of
restricted highest weight in each case. 

For $\lambda = (4,0)$ and $(3,1)$ one could argue as in the preceding
paragraph, or restrict to an appropriate Levi subgroup, as in the last
paragraph of \ref{ss:pf}. 

The remaining cases, up to symmetry, are $\T(3,0)$, $\T(4,1)$, and
$\M$. We apply the translation principle
\cite[II.E.11]{Jantzen}. Observe (from their structure) that $\T(0,2)$
embeds in $\T(4,0)$, which in turn embeds in $\T(2,4)$. Picking
$\lambda = (0,0)$ and $\mu=(-1,1)$ in the closure of the bottom
alcove, observe that applying the (exact) functor $T_\mu^\lambda$ to
these embeddings, we obtain embeddings of $\T(0,3)$ in $\T(4,1)$, and
$\T(4,1)$ in $\T(3,3)$.  Since $\T(3,3)$ has simple $G_1T$-socle of
restricted highest weight, it follows that the same holds for
$\T(0,3)$ and $\T(4,1)$. The cases $\T(3,0)$ and $\T(1,4)$ are treated
by the symmetric argument. Finally, we observe that $\dim_K
\Hom_{G_1T} (\L(0,0), \L(1,1)\otimes \L(1,1)) = 1$, by a calculation
similar to \ref{M}(1). This, along with \ref{M}, shows that $\M$
remains indecomposable on restriction to $G_1T$, with socle and head
isomorphic to $\L(1,1)$, and with 7 copies of $\L(0,0)$ in the middle
Loewy layer.  The proof of Proposition \ref{main:3} is complete.

\subsection{Discussion}\label{disc}\noindent
We now discuss the remaining issue in characteristic 3: the structure
of the tilting modules $\T(\lambda)$ for $\lambda = (3,3)$, $(4,3)$,
$(3,4)$, and $(4,4)$. These tilting modules are in fact $S$-modules
for the Schur algebra $S = S_K(3, r)$ in degree $r = 9, 10, 11, 12$,
respectively. (See \cite{Green, Martin} for background on Schur
algebras.) 

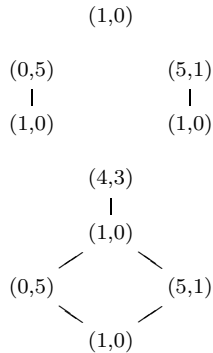
\begin{figure}[ht]
\[
\def\objectstyle{\scriptstyle}
\xymatrix@=6pt{
  & (1,0) & \\
(0,5)\ar@{-}[d] & & (5,1)\ar@{-}[d] \\
(1,0) & & (1,0) \\
  & (4,3) \ar@{-}[d] & \\
  & (1,0) \ar@{-}[dl] \ar@{-}[dr]  & \\
(0,5) \ar@{-}[dr] &  & (5,1) \ar@{-}[dl]\\
 & (1,0) &
}
\]
\caption{\small Weyl filtration factors of $\T(4,3)$}\label{fig:two}
\end{figure}

Thus, in order to study the structure of $\T(\lambda)$ one may employ
techniques from the theory of finite dimensional quasi-hereditary
algebras. Now the simplest cases (in terms of number of composition
factors) are $\T(4,3)$ for $S(3,10)$ and $\T(3,4)$ for $S(3,11)$. As
these modules are symmetric, it makes sense to focus on the smaller
Schur algebra $S(3,10)$ and thus $\T(4,3)$. In fact, it is enough to
understand the block $A$ of $S(3,10)$ consisting of the six weights
$(10,0)$, $(6,2)$, $(4,3)$, $(5,1)$, $(0,5)$, and $(1,0)$. (It is
easily seen that this is a complete linkage class of dominant weights
in $S(3,10)$, for instance by drawing the alcove diagrams.)  To
construct $\T(4,3)$ we must ``glue'' together the $\Delta$-factors in
a way that results in a contravariantly self-dual module. Looking at
the diagrams in Figure \ref{fig:two} below
picturing the various Weyl modules in the filtration, we see that it
is impossible to do this in a rigid way. There are three copies of
$\L(1,0)$ above the middle factor $\L(4,3)$ and only two below. Thus,
there must be two copies of $\L(1,0)$ lying immediately above
$\L(4,3)$ when viewing the radical series, and two copies lying
immediately below $\L(4,3)$ when viewing the socle series. This
implies that $\T(4,3)$ is not rigid.  To understand the structure of
$\T(4,3)$ one may apply Gabriel's theorem to find a quiver and
relations presentation for the basic algebra of the block $A$, or an
appropriate quasi-hereditary quotient thereof. This is carried out in
the Appendix. The other cases could be treated similarly.

Note that none of $\T(4,3)$, $\T(3,4)$, $\T(3,3)$, or $\T(4,4)$ is
projective as an $S$-module, because if so, the reciprocity law
$(P(\lambda)\colon \Delta(\mu)) = [\nabla(\mu) \colon L(\lambda)]$
(see e.g.\ \cite[Prop.~A2.2]{Donkin:book}) would be violated.


%
%



\newcommand{\rad}{\mbox{\rm rad}}

\newcommand{\topp}{\mbox{\rm top}}

\newcommand{\Ker}{\mbox{\rm Ker}}

\newcommand{\Image}{\mbox{\rm Im}}

\renewcommand{\mod}{\mbox{\rm mod}}

\def\arr#1#2{\arrow <1.5mm> [0.25,0.75] from #1 to #2}




%




\newpage
\appendix
\thispagestyle{empty}

\vspace{1in}{\ }

\centerline{\large\bf APPENDIX: THE $\SL_3$-MODULE $T(4,3)$ FOR $p=3$}
\smallskip

\centerline{\large C.~M.~Ringel}

\bigskip

\noindent
Let $k$ be an algebraically closed field of characteristic
$p=3$. Following Bowman, Doty and Martin, we consider rational
$\SL_3$-modules with composition factors $L(\lambda)$, where $\lambda$
is one of the weights $(1,0),\; (0,5),\; (5,1),\; (4,3),\; (6,2)$.
Dealing with a dominant weight $(a,b)$, or the simple module $L(a,b)$,
we usually will write just $ab$.  The corresponding Weyl module, dual
Weyl module, or tilting module, will be denoted by $\Delta(ab),\;
\nabla(ab)$ and $T(ab)$, respectively.

The paper [BDM] by Bowman, Doty and Martin describes in detail the
structure of the modules $\Delta(\lambda), \nabla(\lambda)$ for
$\lambda = 10,\; 05,\; 51,\; 43,\; 62$ and also $T(10),\; T(05),\;
T(51)$ and it provides the factors of a $\Delta$-filtration for
$T(43)$. This module $T(43)$ is still quite small (it has length
$10$), but its structure is not completely obvious at first sight. The
main aim of this appendix is to explain the shape of this module.

Let us call a finite set $I$ of dominant weights (or of simple
modules) an {\it ideal} provided for any $\lambda \in I$ all
composition factors of $T(\lambda)$ belong to $I$. The category of
modules with all composition factors in an ideal $I$ is a highest
weight category with weight set $I,$ thus can be identified with the
module category of a basic quasi-hereditary algebra which we denote by
$A(I)$.  In order to analyse the module $T(43)$, we need to look at
the ideal $I = \{10,05,51,43\},$ thus at the algebra $A(10,05,51,43)$.

In order to determine the precise relations for $A(10,05,51,43),$
we will have to look also at the module $T(62)$, see section 4.
Note that $\{ 10,\; 05,\; 51,\; 43,\; 62\}$
is again an ideal, thus we deal with the
algebra $A(10,05,51,43,62).$

The use of quivers and relations for presenting a basic finite
dimensional algebras was initiated by Gabriel around 1970, 
the text books [ARS] and [ASS] can be used as a reference. 
The class of quasi-hereditary algebras
was introduced by Scott and Cline-Parshall-Scott; for basic properties
one may refer to [DR] and [R2].
The author is grateful to S. Doty and R. Farnsteiner for fruitful 
discussions and helpful suggestions
concerning the material presented in the appendix.

\pagestyle{myheadings}
\markboth{\thepage\hfill C.M.~RINGEL\hfill}%
{\hfill APPENDIX: THE $\SL_3$-MODULE $T(4,3)$ FOR $p=3$\hfill\thepage}
\section*{1. The main result}

\noindent
Deviating  from [BDM],  we will  consider {\bf  right}  modules.  Thus,
given a  finite-dimensional algebra $A$,  an indecomposable projective
$A$-module is of  the form $eA$ with $e$  a primitive idempotent.  The
algebras to be considered will  be factor algebras of path algebras of
quivers and the advantage of looking  at right modules will be that in
this way we  can write the paths  in the quiver as going  from left to
right.
       \medskip

{\bf Proposition.} {\it The algebra $A(10,05,51,43)$ is isomorphic to
  the path algebra of the quiver
$$
\hbox{\beginpicture
\setcoordinatesystem units <1cm,1cm>
\put{$Q = Q(10,05,51,43)$} at -3 0
\put{$10$} at 0 0
\put{$05$} at 1 1 
\put{$51$} at 1 -1 
\put{$43$} at 2 0 
\arr{0.3 0.08}{1.7 0.08}
\arr{1.7 -.08}{0.3 -.08}
\arr{0.1 0.3}{0.7 0.8}
\arr{0.8 0.7}{0.2 0.2}
\arr{0.7 -.8}{0.1 -.3}
\arr{0.2 -.2}{0.8 -.7}
\put{$\alpha$} at 0.2 0.8
\put{$\alpha'$} at 0.8 0.5
\put{$\beta'$} at 0.2 -.8
\put{$\beta$} at 0.8 -0.5
\put{$\gamma$} at 1.3 0.3
\put{$\gamma'$} at 1.3 -0.3
\endpicture}
$$
modulo the ideal generated by the following
relations}
$$
\begin{array}{c}
 \alpha'\alpha  = 0, \quad
 \alpha'\beta  = 0, \quad
 \beta'\alpha  = 0, \quad
 \beta'(1-\gamma\gamma')\beta = 0, \\
  \gamma'\gamma  = 0, \quad 
 \gamma'(\alpha\alpha' - \beta\beta') = 0, \quad
 (\alpha\alpha' - \beta\beta')\gamma  = 0, \quad \gamma'\alpha\alpha'\gamma = 0.
\end{array}
$$
	\medskip

We are going to give some comments before embarking on the proof. 
   \medskip

(1) Since the quiver $Q(10,05,51,43)$ 
is bipartite, say with a $(+)$-vertex $10$ and three $(-)$-vertices
$05,\; 51,\; 43$, possible relations between vertices of the same parity involve paths of
even lengths, those between vertices with different parity involve paths of odd lengths.
Our convention for labelling arrows between a $(+)$-vertex $a$ 
and a $(-)$-vertex $b$ is the following: we use a greek letter for the arrow 
$a \to b$ and add a dash for the arrow $b \to a$.
   \medskip

(2) The assertion of the proposition  can be visualized by drawing the 
shape of the indecomposable projective $A$-modules.
The  indecomposable projective
$A$-module with top $\lambda$ will be denoted by $P(\lambda) = e_\lambda A$,
where $e_\lambda$ is the primitive idempotent corresponding to $\lambda$,
and we will denote the radical of $A$ by $J$.
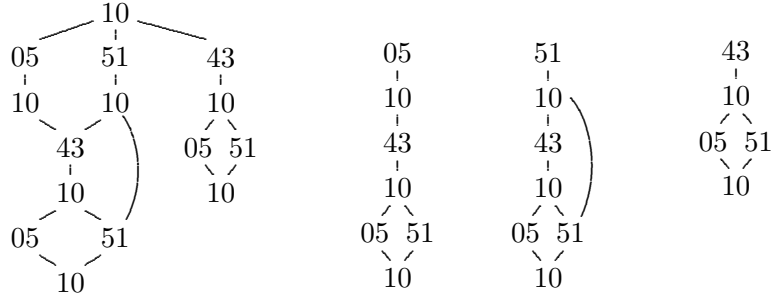
$$
\hbox{\beginpicture
\setcoordinatesystem units <1cm,.6cm>
\put{\beginpicture

\plot 0  .6  0  .4 /
\plot 1.2 .6  1.2 .4 /

\put{$05$} at 0 1
\put{$10$} at 0 0
\put{$51$} at 1.2 1
\put{$10$} at 1.2 0

\plot 0.4 -0.6  0.2 -0.4 /
\plot 0.8 -0.6  1 -0.4 /

\put{$43$} at 0.6 -1
\plot 0.6 -1.4  0.6 -1.6 /

\put{$10$} at 0.6 -2
\plot 0.4 -2.4  0.2 -2.6 /
\plot 0.8 -2.4  1 -2.6 /
\plot 0.4 -3.6  0.2 -3.4 /
\plot 0.8 -3.6  1 -3.4 /
\ellipticalarc axes ratio 1:1.8 -100 degrees from 1.3 -.3 center at 1 -1.5

\put{$05$} at 0   -3
\put{$51$} at 1.2 -3
\put{$10$} at 0.6 -4

\put{$43$} at 2.6 1
\put{$10$} at 2.6 0
\put{$05$} at 2.3 -1
\put{$51$} at 2.9 -1
\put{$10$} at 2.6 -2
\plot 2.6 .6  2.6 .4 /
\plot 2.5 -.4    2.4 -.6 /
\plot 2.7 -.4    2.8 -.6 /
\plot 2.5 -1.6    2.4 -1.4 /
\plot 2.7 -1.6    2.8 -1.4 /

\plot 1.2 1.6  1.2 1.4 /
\plot 0.9 1.8  0.2 1.4 /
\plot 1.5 1.8   2.4 1.4 /
\put{$10$} at 1.2 2

\endpicture} at -1.5 0 

\put{\beginpicture
\plot -0.1 -1.4  -0.2 -1.6 /
\plot   .1 -1.4    .2 -1.6 /
\put{} at 0 3
\put{$05$} at 0 2
\put{$10$} at 0 1
\put{$43$} at 0 0
\put{$10$} at 0 -1
\put{$05$} at -.3 -2
\put{$51$} at .3 -2
\plot 0 1.4  0 1.6 /
\plot 0 0.4  0 0.6 /
\plot 0 -.4  0 -.6 /

\plot -0.1 -2.6  -0.2 -2.4 /
\plot   .1 -2.6    .2 -2.4 /
\put{$10$} at 0 -3 
\endpicture} at 2 0 

\put{\beginpicture
\plot -0.1 -1.4  -0.2 -1.6 /
\plot   .1 -1.4    .2 -1.6 /
\put{} at 0 3
\put{$51$} at 0 2
\put{$10$} at 0 1
\put{$43$} at 0 0
\put{$10$} at 0 -1
\put{$05$} at -.3 -2
\put{$51$} at .3 -2
\plot 0 1.4  0 1.6 /
\plot 0 0.4  0 0.6 /
\plot 0 -.4  0 -.6 /

\ellipticalarc axes ratio 1:1.8 -100 degrees from 0.3 1 center at 0 -.5
\plot -0.1 -2.6  -0.2 -2.4 /
\plot   .1 -2.6    .2 -2.4 /
\put{$10$} at 0 -3 

\endpicture} at 4 0 

\put{\beginpicture
\put{} at 0 -2.5

\put{$43$} at 0 2
\put{$10$} at 0 1
\put{$05$} at -0.3 0
\put{$51$} at 0.3 0
\put{$10$} at 0 -1
\plot -0.1 0.6  -0.2 0.4 /
\plot   .1 0.6    .2 0.4 /
\plot -0.1 -.6  -0.2 -.4 /
\plot   .1 -.6    .2 -.4 /
\plot 0 1.4  0 1.6 /

\endpicture} at 6.5 0 

\endpicture}
$$ These are the coefficient quivers of the indecomposable projective
$A$-modules with respect to suitable bases. In addition, the
proposition asserts that all the non-zero coefficents can be chosen to
be equal to $1$.  Note that this means that $A$ has a basis $\mathcal
B $ which consists of a complete set of primitive and orthogonal
idempotents as well as of elements from the radical $J$, and such that
$\mathcal B$ is multiplicative (this means: if $u,v\in \mathcal B,$
then either $uv = 0$ or else $uv \in \mathcal B $).

For the convenience
of the reader, let us recall the notion of a
coefficient quiver (see for example [R3]): 
By definition, a representation $M$ of a quiver $Q$  over a field $k$ 
is of the form
$M = (M_x;M_\alpha)_{x,\alpha}$; here, for every vertex $x$ of $Q$, 
there is given
a finite-dimensional $k$-space $M_x$, say of dimension $d_x$, and
for every arrow $\alpha: x \to y,$
there is given a linear transformation $M_\alpha: M_x \to M_y$. 
A {\it basis\/} $\mathcal B$ of $M$ is by definition a subset of the 
disjoint union 
of the various $k$-spaces $M_x$ such that for any vertex $x$ the set 
$\mathcal B_x = \mathcal B\cap M_x$ is a basis of $M_x.$ 
Now assume that there is given a basis $\mathcal B$
of $M$. For any arrow $\alpha:x \to y,$ 
write $M_\alpha$ as a $(d_x\times d_y)$-matrix $M_{\alpha,\mathcal B}$ 
whose rows are indexed by $\mathcal B_x$
and whose columns are indexed by $\mathcal B_y.$ We denote by 
$M_{\alpha,\mathcal B}(b,b')$ 
the corresponding matrix coefficients, where $b\in \mathcal B_x,$
$b'\in \mathcal B_y,$ 
these matrix coefficients  $M_{\alpha,\mathcal B}(b,b')$ are defined by
$M_\alpha(b) = \sum_{b'\in \mathcal B} b'\,M_{\alpha,\mathcal B}(b,b').$
By definition, the {\it coefficient quiver\/} $\Gamma(M,\mathcal B)$ 
of $M$ with respect to $\mathcal B$ has 
the set $\mathcal B$ as set of vertices, and there is an arrow $(\alpha,b,b')$
provided $M_{\alpha,\mathcal B}(b,b') \neq 0$ (and we call 
$M_{\alpha,\mathcal B}(b,b')$ the corresponding coefficient).
If $b$ belongs to $\mathcal B_x$,
we usually label the vertex $b_x$ by $x$. If necessary, we label
the arrow $(\alpha,b,b')$ by $\alpha$; but since we only deal with quivers without
multiple arrows, the labelling of arrows could be omitted. 
In all cases considered in the
appendix, we can arrange the vertices in such a way that all the arrows
point downwards, and then replace arrows by edges. 
This convention will be
used throughout.

Note that there is a
long-standing tradition in matrix theory to focus attention to such 
coefficient quivers (see e.g. [BR]), whereas the representation theory of
groups and algebras is quite reluctant to use them.

Looking at the pictures one should be aware that the four upper
base elements form a complete
set of primitive and orthogonal idempotents, thus these 
are the generators of the indecomposable projective $A$-modules.
Those directly below generate the radical of $A$, and they are
just the arrows of the quiver (or better: the residue classes of the arrows
in the factor algebra of the path algebra modulo the relations). 
Of course, on the left we see $P(10)$, then $P(05)$ and $P(51)$,
and finally, on the right, $P(43)$.
    \medskip

(3) The strange relation $\beta'(1-\gamma\gamma')\beta = 0$ leads to the
curved edge in $P(51)$ as well as in $P(10)$. 
Note that the submodule lattice of $P(51)$
would not at all be changed when deleting this extra line --- but its effect 
would be seen in $P(10)$. Namely, 
without this extra line, the socle of $P(10)$ would be of length 
$3$ (namely, $\topp\; \rad^2 P(10)$ is the 
direct sum of three copies of $10$, and the two copies displayed in the left part
are both mapped under $\gamma$ to $43$, thus there is a diagonal which
is mapped under $\gamma$ to zero; without the curved line, this diagonal would belong
to the socle), whereas the socle of $P(10)$ is of length $2$.
   \medskip

(4) Looking at the first four
relations presented above, one could have the feeling of 
a certain asymmetry concerning the role of $P(05)$ and $P(51)$,
or also of the role of $05$ and $51$ 
as composition factors of the radical of $P(51).$
But such a feeling is misleading as will be seen in the proof. 
The pretended lack of symmetry concerns also our display of $T(43)$.
Sections 7 and 8 will be devoted to a detailed analysis of the module $T(43)$ 
in order to focus the attention to its hidden symmetries.
   \medskip

(5) Note that all the tilting $A$-modules are local (and also colocal):
$$
\begin{array}{ccccc}
 T(10) &=& P(10)/(\alpha A+\beta A+ \gamma A) \\
 T(05) &=& P(10)/(\beta A+\gamma A),\quad\\ 
 T(51) &=& P(10)/(\alpha A+\gamma A),\quad\\
T(43) &=& P(10)/\gamma A. \qquad
\end{array}
$$
As we have mentioned, sections 7 and 8
will discus in more detail the module $T(43).$
     \medskip

(6) A further comment: One may be surprised to see that one can find relations	
which are not complicated at all: many are monomials, the remaining ones are 
differences of monomials, always using paths of length at most 4.

\section*{2. Preliminaries on algebras and the presentation of algebras 
      using quivers and relations}

\noindent 
Let $t$ be a natural number.  Recall that the zero module has Loewy
length $0$ and that a module $M$ is said to have {\it Loewy length at
  most $t$} with $t\ge 1$, provided it has a submodule $M'$ of Loewy
length at most $t-1$ such that $M/M'$ is semisimple.  Given a module
$M$, we denote by $\soc_tM$ the maximal submodule of Loewy length at
most $t$, and by $\topp^tM$ the maximal factor module of Loewy length
$t$. Of course, we write $\soc = \soc_1$ and $\topp = \topp^1,$ but
also $\topp^t M = M/\rad^t M.$

Let $A$ be a finite-dimensional basic algebra with radical $J$ and quiver $Q$.
Let us
assume that $Q$ has no multiple arrows (which is the case for all
the quivers considered here). 
For any arrow $\zeta:i \to j$ in $Q$, we choose an element 
$\eta(\zeta)\in e_iJe_j\setminus e_iJ^2e_j$;
the set of elements $\eta(\zeta)$ will be called a {\it generator choice} for $A$.
In this way, we obtain a surjective algebra homomorphisms 
$$
 \eta:kQ \to A
$$ 
If $\rho$ is the kernel of $\eta$, then $\rho = \bigoplus_{ij} e_i\rho e_j$,
and we call a generating set for $\rho$ consisting of elements in 
$\bigcup_{ij} e_i\rho e_j$ a {\it set of relations} for $A$. We are looking
for a generator choice for the algebra
$A(10,05,51,43)$ which allows to see clearly the structure of $T(43)$.
Usually, we will write $\zeta$ instead of $\eta(\zeta)$ and 
hope this will not produce confusion. 
If $\zeta\in e_iJe_j\setminus e_iJ^2e_j$ belongs to a generator choice, we
obviously may replace it by any element of the form $c\zeta + d$ with $0\neq c\in
k$ and $d\in e_iJ^2e_j$ and obtain a new generator choice.

\section*{3. The algebra $B = A(10,05,51)$}

\noindent
Consider a quasi-hereditary algebra
$B$ with quiver being the full subquiver of $Q(10,05,51,43)$ with vertices
$10,\; 05,\; 51$ and with ordering $10 < 05,$
$10 < 51$. It is well-known (and easy to see) that $B$ is uniquely determined
by these data. The indecomposable projectives have the following shape
$$
\hbox{\beginpicture
\setcoordinatesystem units <1cm,.6cm>

\put{\beginpicture
\put{$05$} at 0 2
\put{$10$} at 0 1
\plot 0 1.4  0 1.6 /  
\endpicture} at 2 0 

\put{\beginpicture
\put{$51$} at 0 2
\put{$10$} at 0 1
\plot 0 1.4  0 1.6 /  
\endpicture} at 4 0 

\put{\beginpicture
\plot 0.3 1.4  0.5 1.7 /
\plot 1.1 1.4  0.9 1.7 /
\put{$10$} at 0.7 2
\put{$05$} at 0 1
\put{$10$} at 0 0
\put{$51$} at 1.4 1
\put{$10$} at 1.4 0
\plot 0 .4  0 .6 /  
\plot 1.4 .4  1.4 .6 /  

\endpicture} at 7 0 

\endpicture}
$$
What we display are the again coefficient quivers of the indecomposable projective
$B$-modules considered as  representations of $kQ$ with respect
to a suitable basis. 
   \bigskip

We see that the algebra $B$ is of Loewy length $3$ and
that it can be described by the relations:
$$
 \alpha'\alpha =  \alpha'\beta = \beta'\alpha =  \beta'\beta = 0.
$$
	\medskip

Of course, $\Delta(10) = \nabla(10) = 10$; and the modules $\Delta(05),$
$\Delta(51),$ $\nabla(05)$ and $\nabla(51)$ are serial of length 2, always
with $10$ as one of the composition factors. This means that
the structure of the modules $\Delta(\lambda),$ $\nabla(\lambda)$,
for $\lambda = 10,\; 05\; 51$ can be read off from the quiver 
(but, of course, conversely, the quiver was obtained from the knowledge of 
the corresponding $\Delta$- and $\nabla$-modules). 

Note that $T(05)$ is the only indecomposable module with a $\Delta$-filtration
with factors $\Delta(10)$ and $\Delta(05)$, since $\Ext^1(\Delta(10),\Delta(05)) = k.$
Similarly, 
$T(51)$ is the only indecomposable module with a $\Delta$-filtration
with factors $\Delta(10)$ and $\Delta(51)$.

Let us remark that the structure of the module category $\mod\; B$ is
well-known: using covering theory, one observes that $\mod\; B$
is  obtained from the category of representations of the affine quiver
of type $\widetilde A_{22}$ with a unique sink and a unique source by
identifying the simple projective module with the simple injective module.
In $\mod\; B$, there is a family of homogeneous tubes indexed by $k\setminus\{0\}$,
the modules on the boundary are of length $4$ with top and socle equal to
$10$ and with $\rad/\soc = 05 \oplus 51$. We will call these modules the 
{\it homogeneous $B$-modules of length $4$.}
(The representation theory of affine quivers can be
found for example in [R1] and [SS]; from covering theory, 
we need only the process of removing a node, see [M].)

\section*{4. The modules $\rad\; \Delta(43)$ and 
         $\nabla(43)/\soc$ are isomorphic}

\noindent
We will use the following information concerning the modules
$\Delta(43)$ and $\nabla(43)$, see [BDM].  Both $\rad\; \Delta(43)$ and
$\nabla(43)/\soc$ are homogeneous $B$-modules of length $4$, thus the
modules $\Delta(43)$ and $\nabla(43)$ have the following shape
$$
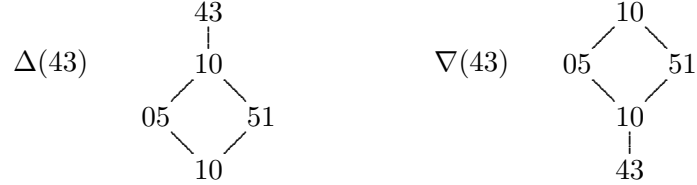
\hbox{\beginpicture
\setcoordinatesystem units <.7cm,.7cm>
\put{\beginpicture
\put{$\Delta(43)$} at -3 2
\put{$10$} at 0 0
\put{$05$} at -1 1
\put{$51$} at 1 1
\put{$10$} at 0 2
\put{$43$} at 0 3
\plot 0 2.3  0 2.7 /
\plot -0.3 0.3  -0.7 0.7 /
\plot -0.7 1.3  -0.3 1.7 /
\plot  0.7 1.3   0.3 1.7 /
\plot  0.3 0.3   0.7 0.7 /
\endpicture} at 0 0 
\put{\beginpicture
\put{$\nabla(43)$} at -3 1

\put{$10$} at 0 0
\put{$05$} at -1 1
\put{$51$} at 1 1
\put{$10$} at 0 2
\put{$43$} at 0 -1
\plot 0 -.3  0 -.7 /
\plot -0.3 0.3  -0.7 0.7 /
\plot -0.7 1.3  -0.3 1.7 /
\plot  0.7 1.3   0.3 1.7 /
\plot  0.3 0.3   0.7 0.7 /

\endpicture} at 8 0 
\endpicture}
$$
Here, we have drawn again coefficient quivers with respect to suitable
bases. But note
that we do not (yet) claim that all the non-zero coefficients can be chosen to
be equal to $1$.

In order to show the assertion in the title, we have to expand our considerations
taking into account also the weight $62$. The existence of an isomorphism 
in question 
will be obtained by looking at the tilting module $T(62).$

In dealing with a tilting module $T(\mu)$, there is a unique submodule 
isomorphic to $\Delta(\mu)$, and a unique factor module isomorphic to $\nabla(\mu).$
Let $R(\mu) = \rad\;\Delta(\mu)$ and let $Q(\mu)$ be the kernel of the canonical map 
$\pi:T(\mu) \to \nabla(\mu)/\soc.$ Note that $\Delta(\mu) \subseteq Q(\mu)$
(namely, if $\pi(\Delta(\mu))$ would not be zero, then it would be a 
submodule of $\nabla(\mu)/\soc$ with top equal to $\mu$; however $\nabla(\mu)/\soc$
has no composition factor of the form $\mu$). It follows that
$R(\mu) \subset Q(\mu)$ and we call $C(\mu) = Q(\mu)/R(\mu)$ the {\it core}
of the tilting module $T(\mu)$. Also, we see that  $\mu = \Delta(\mu)/R(\mu)$  is a
simple submodule of $C(\mu).$ In fact, {\it $\mu$ is a direct summand of $C(\mu)$.}
Namely, there is $U\subset T(\mu)$ with $T(\mu)/U = \nabla(\mu).$ Then 
$U \subset Q(\mu)$ and $Q(\mu)/U = \mu.$ Since $R(\mu) \subset Q(\mu)$
and $R(\mu)$ has no composition factor of the form $\mu$, it follows that
$R(\mu) \subseteq U.$ Altogether, we see that $U+\Delta(\mu) = Q(\mu)$
and $U\cap \Delta(\mu) = R(\mu).$ Thus $Q(\mu)/R(\mu) = U/R(\mu) \oplus 
\Delta(\mu)/R(\mu) = U/R(\mu)\oplus \mu.$ 

The module $\Delta(62)$ is serial with going down factors $62,\; 43,\;
10,\; 51$, and the module $\nabla(62)$ is serial with going down
factors $51,\; 10,\; 43,\; 62$, see [BDM], 4.1. Also we will use
that $T(62)$ has $\Delta$-factors $\Delta(51),\; \Delta(43),\;
\Delta(62),$ each with multiplicity one (and thus $\nabla$-factors
$\nabla(62),\; \nabla(43),\;\nabla(51)$).  To get the $\Delta$-factors
of $T(62)$, one has to use [BDM], (2.2.2) along with the known
structure of the Deltas (this requires a small calculation, which is
left to the reader.)

The quiver $Q(10,05,51,43,62)$ of $A(10,05,51,43,62)$ is
$$
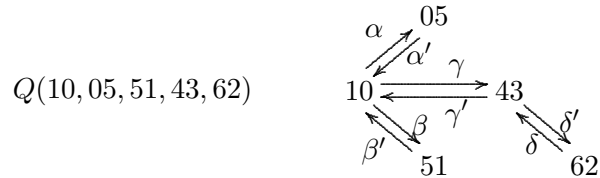
\hbox{\beginpicture
\setcoordinatesystem units <1cm,1cm>
\put{$10$} at 0 0
\put{$05$} at 1 1 
\put{$51$} at 1 -1 
\put{$43$} at 2 0 
\arr{0.3 0.08}{1.7 0.08}
\arr{1.7 -.08}{0.3 -.08}
\arr{0.1 0.3}{0.7 0.8}
\arr{0.8 0.7}{0.2 0.2}
\arr{0.7 -.8}{0.1 -.3}
\arr{0.2 -.2}{0.8 -.7}
\put{$\alpha$} at 0.2 0.8
\put{$\alpha'$} at 0.8 0.5
\put{$\beta'$} at 0.2 -.8
\put{$\beta$} at 0.8 -0.5
\put{$\gamma$} at 1.3 0.3
\put{$\gamma'$} at 1.3 -0.3

\put{$62$} at 3 -1 
\arr{2.7 -.8}{2.1 -.3}
\arr{2.2 -.2}{2.8 -.7}
\put{$\delta$} at 2.3 -.7
\put{$\delta'$} at 2.8 -0.4
\put{$Q(10,05,51,43,62)$} at -3 0
\endpicture}
$$
with ordering $10 < 05 < 43 < 62$, and $10 < 51 < 43$.
     \bigskip

{\bf Lemma 1.} {\it The core of $T(62)$ is of the form $(\rad\;\Delta(43))\oplus 62$
as well as of the form $(\nabla(43)/\soc) \oplus 62.$}
   \medskip

{\bf Corollary.} {\it The modules $\rad\;\Delta(43)$ and $\nabla(43)/\soc$ are isomorphic.}
     \medskip

Note that it is quite unusual that the modules 
$\rad\;\Delta(\lambda)$ and $\nabla(\lambda)/\soc$ are isomorphic, for a weight
$\lambda$. 
	   \medskip

{\em Proof of Lemma 1.} 
Let $T_1 \subset T_2 \subset T(62)$ be a filtration with factors
$$
 T_1 = \Delta(62),\quad T_2/T_1 = \Delta(43),\quad T(62)/T_2 = \Delta(51).
$$
Now $R(62) = \rad\;\Delta(62) \subset T_1 \subset T_2,$ thus we may look at 
the factor module $T_2/R(62)$ and the exact sequence
$$
 0 \to 62 \to T_2/R(62) \to \Delta(43) \to 0
$$
(with $62 = T_1/R(62)$). We consider the submodule 
$N = \rad\;\Delta(43)$ of $\Delta(43)$, with factor module $\Delta(43)/N = 43$. 
 We have
$\Ext^1(N,62) = 0$, since $\Ext^1(S,62) = 0$ for all the
composition factors $S$ of $N.$ This implies that
there is an exact sequence
$$
 0 \to N \oplus 62 \to T_2/R(62) \to 43 \to 0.
$$
Thus, there is a submodule $U \subset T_2$ with $R(62) \subset U$ such that
$U/R(62)$ is isomorphic to $N \oplus 62$ and $T_2/U$ is isomorphic to $43$.
Since $T(62)/T_2 = \Delta(51)$ is of length $2$, we see that $T(62)/U$
is of length $3.$   

Now consider the canonical map $\pi:T(62) \to \nabla(62)/\soc$.
This map vanishes on $R(62)$, thus induces a map 
$\pi': T(62)/R(62) \to \nabla(62)/\soc.$
Let us look at the submodule $U/R(62)$ of $T(62)/R(62)$.
Since the socle
of $\nabla(62)/\soc$ is equal to $43$, and $U/R(62) = N\oplus 62$ has
no composition factor of the form $43$, we see that $U/R(62)$ is contained
in the kernel of $\pi'$, and therefore $U$ is contained in the kernel of $\pi$.

By definition, the kernel of the canonical map $\pi:T(62) \to \nabla(62)/\soc$ is $Q(62)$,
thus we have shown that $U \subseteq Q(62)$.
But $T(62)/U$ is of length 3 as is $T(62)/Q(62),$ thus
$U = Q(62)$. But this means that $Q(62)/R(62) = U/R(62) = N \oplus 62 =  
(\rad\;\Delta(43)) \oplus 62.$

The dual arguments show that $Q(62)/R(62) = (\nabla(43)/\soc) \oplus 62$.
\hfill$\square$
    \medskip

As we have mentioned, 
the module $N = \rad\;\Delta(43)$ is a $B$-module, where $B = A(10,05,51)$.
This algebra $B$ has been discussed in section 3. The coefficient quiver of $N$ is
$$
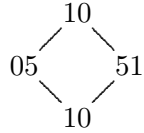
\hbox{\beginpicture
\setcoordinatesystem units <.7cm,.7cm>
\put{$10$} at 0 0
\put{$05$} at -1 1
\put{$51$} at 1 1
\put{$10$} at 0 2
\plot -0.3 0.3  -0.7 0.7 /
\plot -0.7 1.3  -0.3 1.7 /
\plot  0.7 1.3   0.3 1.7 /
\plot  0.3 0.3   0.7 0.7 /
\endpicture}
$$
Now, choosing a suitable basis of $N$, we can assume that at least 3 
of the non-zero coefficients are equal to $1$ and we look at the
remaining coefficient, say that for the arrow $\alpha$. 
It will be a non-zero scalar $c$ in $k$. Recall that we have started with
a particular generator choice for the algebra $B$ which we can change.
If we replace the element $\alpha\in J$ by $\frac1c\alpha$, then 
the coefficients needed for $N$ will all be equal to $1$. 
    \bigskip

{\bf Remark.} Extending the analysis of the $\Delta$- and the
$\nabla$-filtrations of $T(43)$, one can show that 
$T(62)$ is the indecomposable projective $A(10,05,51,43,62)$-module with
top $51$ (as well as the indecomposable injective 
$A(10,05,51,43,62)$-module with
socle $51$). 
As Doty has pointed out, the last assertion follows also from 
Theorem 5.1 of the DeVisscher-Donkin
paper [DD] (that result is
based on their Conjecture 5.2 holding, but it is proved in Section 7
of the same paper that the conjecture holds for $\GL(3);$ hence it holds
also for $\SL(3)$). 

Let us add without proof that in this way one may show that the module $T(62)$ 
has a coefficient quiver of the form
$$
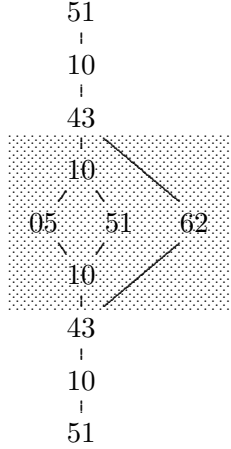
\hbox{\beginpicture
\setcoordinatesystem units <.5cm,.7cm>
\multiput{$10$} at 0 1  0 3  0 5  0 7 /
\multiput{$51$} at 0 0  1 4  0 8 /
\multiput{$43$} at 0 2  0 6 /
\put{$05$} at -1 4 
\put{$62$} at 3 4 
\plot 0 0.4  0 0.6 /
\plot 0 1.4  0 1.6 /
\plot 0 2.4  0 2.6 /
\plot 0 5.4  0 5.6 /
\plot 0 6.4  0 6.6 /
\plot 0 7.4  0 7.6 /
\plot -0.4 3.4  -0.6 3.6 /
\plot -0.4 4.6  -0.6 4.4 /
\plot 0.4 3.4  0.6 3.6 /
\plot 0.4 4.6  0.6 4.4 /
\plot 0.6 2.4  2.6 3.6 /
\plot 0.6 5.6  2.6 4.4 /
\setshadegrid span <.5mm>
\hshade 2.3  -2 4   5.7 -2 4 /
\endpicture}
$$
the shaded part being the core of $T(62)$.

\section*{5. The module $T(43)$}

{\bf Lemma 2.} {\it  We have $\topp\; T(43) = 10 = \soc\; T(43).$}
     \medskip

{\em Proof:} We use that $T(43)$ has $\Delta$-factors 
$\Delta(10),\; \Delta(05),\; \Delta(51),\;
\Delta(43)$ in order to show that $\topp\; T(43) = 10.$ 
Since $\topp\; T(43)$  is isomorphic to
a submodule of the direct sum of the tops of the $\Delta$-factors,
it follows that $\topp\; T(43)$ is multiplicity free.
Since $T(43)$ maps onto $\nabla(43)$,
the only composition factor $43$ cannot belong to the top. 

Actually, it is $N = T(43)/\rad\; \Delta(43)$ which 
maps onto $\nabla(43)$, and $\nabla(43)$ maps onto
$\nabla(05)$ which is serial with top $10$ and socle $05$;
this shows that the only composition factor of the form $05$ of $N$ 
does not belong to $\topp\; N$.
Now $05$ is not in $\topp\; N$ and not in $\topp\;\rad\; \Delta(43)$, thus not in $\topp\; T(43)$.
Similarly, $51$ is not in $\topp\; T(43).$ It follows that $\topp\; T(43) = 10.$

Note that the $\nabla$-factors of $T(43)$ are 
$\nabla(10),\; \nabla(05),\;
\nabla(51),\;\nabla(43).$
Namely, $T(43)$ maps onto $\nabla(43)$, say with kernel $N'$.
The number of composition factors of $N'$ of the form $05,51,10$
is $1,1,3$, respectively. 
Since $N'$ has a $\nabla$-filtration,
its $\nabla$-factors have to be $\nabla(05), \nabla(51)$ and $\nabla(10)$,
each with multiplicity one. In the same way, as we have seen that $T(43)$
has simple top $10$, we now see that it also has simple socle $10$.
\hfill$\square$
    \medskip

Let us add also the following remark:
    \medskip

{\bf Remark.} {\it The module $T(43)$ is a faithful $A$-module}.
     \medskip

{\em Proof:} First of all, we show that the modules 
$T(05)$ and $T(51)$ are both isomorphic to factor modules (and to submodules)
of $T(43)$. 
The $\Delta$-filtration of $T(43)$ shows that $T(43)$ has
a factor module with factors $\Delta(10)$ and $\Delta(05)$. Since this factor module is 
indecomposable,
it follows that it is $T(05)$. Similarly, $T(51)$ is a factor module of $T(43).$
(And dually, $T(05)$ and $T(51)$ are also submodules of $T(43)$). Of course,
also $T(10)$ is a factor module and a submodule of $T(43)$.
It follows that
$T(43)$ is faithful, since the 
direct sum of all tilting modules is always a faithful module 
(it is a ``tilting'' module in the sense used in [R2]).\hfill$\square$


\section*{6. Algebras with quiver $Q(10,05,51,43)$}

\noindent
Let us assume that we deal with a quasi-hereditary algebra $A$ with quiver
$Q(10,05,51,43)$, with ordering $10 < 05 < 43$ and $10 < 51 < 43$ and such that
$\rad\; \Delta(43)$ and $\nabla(43)/\soc$ both are homogeneous $B$-modules
of length $4$.

Since we know the composition factors of all the $A$-modules $\nabla(\lambda)$, 
we can use the reciprocity law in order to see that  
the indecomposable projective modules
have the following $\Delta$-factors (going downwards)
$$
 \begin{array}{ccc} 
         P(43) && \Delta(43) \\
         P(05) && \Delta(05) \mid \Delta(43) \\
         P(51) && \Delta(51) \mid \Delta(43) \\
         P(10) && \Delta(10) \mid \Delta(05)\oplus\Delta(51) \mid 
 \Delta(43)\oplus\Delta(43). 
 \end{array}
$$
We see: 
Since the Loewy length of these factors of $P(10)$ are $1,2,4$, the Loewy length
of $P(10)$ can be at most $7$. Of course, the Loewy length of 
$P(43) = \Delta(43)$ is $4$ and that of $P(05)$ and $P(51)$
is at most $6$. It follows that $J^7 = 0.$ 
   \medskip

Our aim is to contruct
a presentation of $A$ by the quiver $Q$ and suitable relations. As we have mentioned,
for any arrow $\alpha:i \to j$ in $Q$ we choose an element in
$e_iJe_j\setminus e_iJ^2e_j$ which we denote again by $\alpha$, in order to
obtain a surjective algebra homomorphisms 
$$
 \eta:kQ \to A.
$$
	\medskip

Since $J^7 = 0$, we see that {\it all paths of length 
$7$ in the quiver are zero when considered as elements of $A$.}
    \medskip

{\bf Lemma 3.} {\it Any generator choice for $A$ satisfies the conditions
$$
\begin{array}{c}
 \alpha'\alpha, \alpha'\beta,
 \beta'\alpha,
 \beta'\beta \ \in J^4, \cr
 \\
  \gamma'\gamma  = 0, \quad 
 \gamma'(\alpha\alpha' - c_0\beta\beta') = 0, \quad
 (\alpha\alpha' - c_1\beta\beta')\gamma  = 0, \quad \gamma'\alpha\alpha'\gamma = 0.
\end{array}
$$
for some non-zero scalars $c_0,c_1\in k.$}
    \medskip

{\em Proof.} The algebra $B$ considered in section 3
is the factor algebra of $A$ modulo the ideal generated by 
$e_{43}.$ Since we know that the paths
$\alpha'\alpha, \; \alpha'\beta,\; \beta'\alpha,\; \beta'\beta$ are zero in
$B$, they belong to $J^4$ (any path between vertices of the form $05$ and $51$
which goes through $43$ has length at least $4$):
$$
 \alpha'\alpha, \alpha'\beta,
 \beta'\alpha,
 \beta'\beta \ \in J^4.
$$
	\medskip

Since $e_{43}Je_{43} = 0,$ we have 
$$
 \gamma'\gamma = 0.
$$
Also, the shape of $P(43)$ shows that $e_{43}J^3e_{10}$ is one-dimensional,
and that the paths $\gamma'\alpha\alpha'$ and $\gamma'\beta\beta'$ both are non-zero,
thus they are scalar multiples of each other. Thus, we can assume that
$$
 \gamma'(\alpha\alpha' - c_0\beta\beta') = 0,
$$
with some non-zero scalar $c_0$. 
Dually, we have
$$
 (\alpha\alpha' - c_1\beta\beta')\gamma = 0
$$
with some non-zero scalar $c_1$. (Later, we will use the fact that
the modules $\rad\; \Delta(43)$ and $\nabla(43)/\soc$ are isomorphic, then we
can assume that $c_0 = c_1$; also, we will replace one of the arrows 
$\alpha, \alpha', \beta, \beta'$ by a
non-zero scalar multiples, in order to change the coefficient $c_0$ to $1$).
	 
Since $P(43) = \Delta(43)$ is of Loewy length $4$, we see that $\gamma'J^3 = 0$,
in particular we have
$$
 \gamma'\alpha\alpha'\gamma = 0
$$
(and also that $\gamma'\alpha\alpha'\alpha$ and 
$\gamma'\alpha\alpha'\beta$ are zero.) \hfill$\square$
			    \medskip

We have seen in the proof that $\gamma'J^3 = 0$, since $\Delta(43)$ is of Loewy length 
$4$. Dually, since $\nabla(43)$ is of Loewy
length $4$, we have $J^3\gamma = 0.$
       \bigskip

{\bf Lemma 4.} 
{\it A factor algebra of the path algebra of the quiver $Q(10,05,51,43)$
satisfying the relations exhibited in Lemma 3 is generated as a $k$-space by the
elements
$$
 \begin{array}{ll} 
 Q_0 & 10,\ 05,\ 51,\ 43, \\
 Q_1 & \alpha,\ \beta,\ \gamma,\ \alpha',\ \beta',\ \gamma',\\
 Q_2 & \alpha\alpha',\ \beta\beta',\ \gamma\gamma',\ \alpha'\gamma,\ 
       \beta'\gamma,\ \gamma'\alpha,\ \gamma'\beta, \\
 Q_3 & \alpha\alpha'\gamma,\ \gamma\gamma'\alpha,\ \gamma\gamma'\beta,\ 
      \alpha'\gamma\gamma',\ \beta'\gamma\gamma',\ \gamma'\alpha\alpha',\\
 Q_4 & \alpha\alpha'\gamma\gamma',\  \gamma\gamma'\alpha\alpha',\ 
      \alpha'\gamma\gamma'\alpha,\ \alpha'\gamma\gamma'\beta,\ 
     \beta'\gamma\gamma'\alpha,\ \beta'\gamma\gamma'\beta, \\ 
 Q_5 & \alpha\alpha'\gamma\gamma'\alpha,\ \alpha\alpha'\gamma\gamma'\beta,\  
     \alpha'\gamma\gamma'\alpha\alpha',\
     \beta'\gamma\gamma'\beta\beta', \\ 
 Q_6 & \alpha\alpha'\gamma\gamma'\alpha\alpha',
 \end{array}
$$
thus is of dimension at most $34$.}
     \medskip

{\em Proof:} One shows inductively that the elements listed as $Q_i$ generate
the factor space $J^i/J^{i+1}.$ This is obvious for $i = 0,1,2$, since here
we have listed all the paths of length $i$. For $i = 3$, the missing
paths of length $3$ are 
$$
 \alpha\alpha'\alpha,\ \alpha\alpha'\beta,\ 
 \beta\beta'\alpha,\ \beta\beta'\beta,\  \gamma\gamma'\gamma,\ 
$$
as well as
$$
 \beta\beta'\gamma,\ 
 \gamma'\beta\beta'.
$$
By assumption, the first five belong to $J^4$, whereas the last two are
equal to a non-zero multiple of $ \alpha\alpha'\gamma$ and
$\gamma'\alpha\alpha'$, respectively.

Next, consider $i\ge 4.$ 
We have to take the paths in $Q_{i-1}$ and multiply them from the right
by the arrows and see what happens. For $i=4$, the missing paths
are $\gamma\gamma'\beta\beta'$ (it is a multiple of $\gamma\gamma'\alpha\alpha'$),
the paths $ \alpha'\gamma\gamma'\gamma$ and $\beta'\gamma\gamma'\gamma$
(both involve $\gamma'\gamma$) as well as the right multiples of
$\gamma'\alpha\alpha'$ (all belong to $J^5$). 

In the same way, we deal with
the cases $i=5, 6, 7$. In particular, for $i=7$, we see that $J^7 = J^8$,
and therefore $J^7 = 0.$ This shows that we have obtained a 
generating set of the algebra as a $k$-space. \hfill$\square$

\section*{7. The algebra $A = A(10,05,51,43)$}

\noindent
Now, let $A = A(10,05,51,43).$ 
     \medskip

{\bf Lemma 5.} {\it For any generator choice of elements of $A$,
the paths listed in Lemma 4 form a basis of $A$.}
    \medskip

{\em Proof:} Lemma 3 asserts that we can apply Lemma 4. On the other hand, we
know that $\dim A = 34,$ since we know the dimension of the indecomposable
projective $A$-modules. \hfill$\square$
	   \bigskip

{\bf Lemma 6.} {\it The socle of $P(10)$ has length $2$.}
     \medskip

{\em Proof.} Since $\Delta(43)\oplus \Delta(43)$ is a submodule of $P(10)$, the
length of the socle of $P(10)$ is at least $2$. 

According to Lemma 2, the top of $T(43)$ is
equal to $10$, thus we see that $T(43)$ is a factor module of $P(10)$, say
$T(43) = P(10)/W$ for some submodule $W$ of $P(10).$
The subcategory of modules
with a $\Delta$-filtration is closed under kernels of surjective
maps [R2], thus $W$ has a $\Delta$-filtration. But $W$ has a composition
factor of the form $43$, and is of length $5$, thus $W$ is isomorphic to 
$\Delta(43)$ and therefore has simple socle. Quoting again Lemma 2, we know
that also $T(43)$ has simple socle, thus the length of the socle of $P(10)$ 
is at most $2$. \hfill$\square$
   \bigskip

{\bf Proof of the proposition.} 
Assume that there is given a generator choice for $A$. Then
$\alpha'\alpha$ belongs to $J^4$, thus to $e_{05}J^4e_{05}.$
The basis of $A$ exhibited in Lemma 4 shows that  
$e_{05}J^4e_{05}$ is generated by $\alpha'\gamma\gamma'\alpha,$
thus we see that $\alpha'\alpha$ has to be a multiple of 
$\alpha'\gamma\gamma'\alpha$. In the same way, we consider also
the elements 
$ \alpha'\beta,\  \beta'\alpha,\ \beta'\beta$ 
and obtain 
scalars $c_{aa},c_{ab},c_{ba},c_{bb}$ (some could be zero) such that
$$
\begin{array}{ccc}
 \alpha'\alpha &=& c_{aa}\,\alpha'\gamma\gamma'\alpha, \cr
 \alpha'\beta &=& c_{ab}\,\alpha'\gamma\gamma'\beta, \cr
 \beta'\alpha &=& c_{ba}\,\beta'\gamma\gamma'\alpha, \cr
 \beta'\beta &=& c_{bb}\,\beta'\gamma\gamma'\beta.
\end{array}
$$

We show that we can achieve that three of these coefficients are zero:
Let 
$$
\begin{array}{ccc}
 \alpha'_0 &=& \alpha'(1-c_{aa}\gamma\gamma'), \qquad\cr
 \beta'_0 &=& \beta'(1-c_{ba}\gamma\gamma'), \qquad \cr
 \beta_0 &=& (1-(c_{ab}-c_{aa})\gamma\gamma')\beta, 
\end{array}
$$
Then 
$$
\begin{array}{ccc}
 \alpha'_0\alpha &=& \alpha'(1-c_{aa}\gamma\gamma')\alpha = 0,\\
  \beta'_0\alpha &=& \beta'(1-c_{ba}\gamma\gamma')\alpha = 0, 
\end{array}
$$
and 
$$
\begin{array}{ccc}
 \alpha'_0\beta_0 &=& \alpha'(1-c_{aa}\gamma\gamma')(1-(c_{ab}-c_{aa})\gamma\gamma')\beta \cr
                   &=& \alpha'(1- c_{aa}\gamma\gamma'-(c_{ab}-c_{aa})\gamma\gamma')\beta \quad\ \cr
                   &=& \alpha'(1- c_{ab}\gamma\gamma')\beta  \quad = \quad 0. \qquad\qquad\quad
\end{array}
$$
In the last calculation, we have deleted the summand in $\rad^6$, since actually
$\gamma'\gamma = 0.$

This shows that replacing $\alpha',\;\beta,\;\beta'$ by $\alpha'_0,\;\beta_0,\;\beta'_0$,
respectively, we can assume that all the parameters $c_{aa}, c_{ab}, c_{ba}$ are
equal to zero. 
      \bigskip

Thus, we can assume that we deal with the relations:
$$
\begin{array}{c}
 \alpha'\alpha  = 0, \quad
 \alpha'\beta  = 0, \quad
 \beta'\alpha  = 0, \quad
 \beta'(1-c_{bb}\gamma\gamma')\beta = 0, \\
  \gamma'\gamma  = 0, \quad 
 \gamma'(\alpha\alpha' - c_0\beta\beta') = 0, \quad
 (\alpha\alpha' - c_1\beta\beta')\gamma  = 0, \quad \gamma'\alpha\alpha'\gamma = 0.
\end{array}
$$

Let us show that $c_{bb}\neq 0.$ Assume, for the contrary that $c_{bb} = 0$. Then
the element $\alpha\alpha'-\beta\beta'$
belongs to the socle of $P(10)$. But of course, also the elements 
$\alpha\alpha'\gamma\gamma'\alpha\alpha'$
and $\gamma\gamma'\alpha\alpha'$ belong to the socle of $P(10)$, thus
the socle of $P(10)$
is of length at least $3$. But this contradicts Lemma 6.

We have mentioned already, that the isomorphy of $\rad\;\Delta(43)$
and $\nabla(43)/\soc$ implies that $c_0 = c_1.$ Thus we deal with
a set of relations
$$
\begin{array}{c}
 \alpha'\alpha  = 0, \quad
 \alpha'\beta  = 0, \quad
 \beta'\alpha  = 0, \quad
 \beta'(1-c'\gamma\gamma')\beta = 0, \\
  \gamma'\gamma  = 0, \quad 
 \gamma'(\alpha\alpha' - c\beta\beta') = 0, \quad
 (\alpha\alpha' - c\beta\beta')\gamma  = 0, \quad \gamma'\alpha\alpha'\gamma = 0.
\end{array}
$$
with two non-zero scalars $c, c'$. It remains a last
change of the generator choice: Replace say $\gamma$ by $\frac 1{c'}\gamma$
and $\alpha$ by $\frac1c\alpha$. Then we obtain the wanted presentation.
This completes the proof of the Proposition. \hfill$\square$

\section*{8. The module $T(43)$}
\noindent
As we have mentioned, $T(43)$ is a factor module of $P(10)$, namely
$T(43) = P(10)/\gamma A$, thus it has the following
coefficient quiver: 
$$
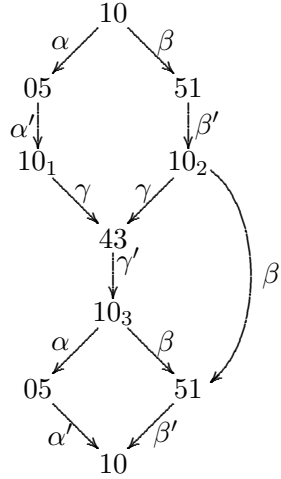
\hbox{\beginpicture
\setcoordinatesystem units <1cm,1cm>
\multiput{$10$} at 0 0   0 6 /
\put{$10_3$} at 0 2 
\put{$10_1$} at -1 4  
\put{$10_2$} at  1 4
\multiput{$05$} at -1 1  -1 5 /
\multiput{$51$} at  1 1   1 5 /
\put{$43$} at 0 3
\arr{-0.8 0.8}{-0.2 0.2}
\arr{-0.2 1.8}{-0.8 1.2}
\arr{0.2 1.8}{0.8 1.2}
\arr{0.8 0.8}{0.2 0.2}
\arr{0 2.8}{0 2.2}

\arr{-0.8 3.8}{-0.2 3.2}
\arr{-0.2 5.8}{-0.8 5.2}
\arr{0.2 5.8}{0.8 5.2}
\arr{0.8 3.8}{0.2 3.2}
\arr{-1 4.8}{-1 4.2}

\arr{1 4.8}{1 4.2}

\ellipticalarc axes ratio 1:1.8 -138 degrees from 1.3 3.9 center at 1 2.5
\arr{1.4 1.18}{1.3 1.08}

\multiput{$\alpha$} at -.7 1.6     -.7 5.6 /
\multiput{$\alpha'$} at -.7 0.4  -1.2 4.5 /
\multiput{$\beta$} at .7 1.6   2.1 2.5     .7 5.6 /
\multiput{$\beta'$} at .7 0.4  1.25 4.5 /
\multiput{$\gamma$} at -0.4 3.6  0.4 3.6 /
\multiput{$\gamma'$} at 0.2 2.7 /

\endpicture}
$$ 
with all non-zero coefficients being equal to $1$.
     \medskip

The picture shows nicely the $\Delta$-filtration of $T(43)$,
but, of course, one also wants to see a $\nabla$-filtration. 
This is the reason why we have labelled the three copies of 
$10$ in the middle (since we exhibit a
coefficient quiver, these elements
$10_1,\ 10_2,\ 10_3$ are elements of a basis).
Consider the subspace
$$
 V = \langle 10_1, 10_2, 10_3 \rangle
$$
of $T(43)$ and the elements $x =  10_1+10_3-10_2$ 
and $y =  10_1-10_2$ of $V$.
One easily sees the following: 

The element $x$ lies in the kernel both of 
$\beta$ and $\gamma$, and it is mapped under $\alpha$ to the 
composition factor $05$ lying in $\soc_2T(43)$. Thus, it provides
an embedding  of $\nabla(05)$ into $T(43)/\soc.$ 

The element $y$ lies in the kernel both of 
$\alpha$ and $\gamma$, and it is mapped under $\beta$ to the 
composition factor $51$ lying in $\soc_2T(43)$.
Thus, it provides
an embedding  of $\nabla(51)$ into $T(43)/\soc.$ 

The sum of the submodules $xA$ and $yA$ is a submodule of $T(43)$ of length
$5$ with a $\nabla$-filtration with factors going down:
$$
 \nabla(05)\oplus \nabla(51) \mid \nabla(10).
$$
Finally, the factor module $T(43)/(xA+yA)$ is obviously of the form
$\nabla(43)$, since its socle is $43$ and its length is $5$.
	      \bigskip

{\bf Remark.} In terms of the
basis of $A$ presented above, we also can write:
$$
\begin{array}{l}
x = \alpha\alpha'+\alpha\alpha'\gamma\gamma'-\beta\beta'\\
y = \alpha\alpha'-\beta\beta'.
\end{array}
$$

\section*{9. A further look at the module $T(43)$}

\noindent
In order to understand the module $T(43)$ better, let us concentrate on the 
essential part which looks quite strange, namely the three subfactors
$10_1,\; 10_2,\; 10_3$ shaded below:
$$
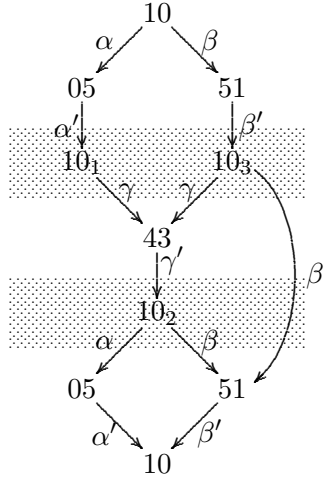
\hbox{\beginpicture
\setcoordinatesystem units <1cm,1cm>
\multiput{$10$} at 0 0  0 6 /

\put{$10_1$} at -1 4  
\put{$10_2$} at 0 2
\put{$10_3$} at  1 4 
\multiput{$05$} at -1 1  -1 5 /
\multiput{$51$} at  1 1   1 5 /
\put{$43$} at 0 3
\arr{-0.8 0.8}{-0.2 0.2}
\arr{-0.2 1.8}{-0.8 1.2}
\arr{0.2 1.8}{0.8 1.2}
\arr{0.8 0.8}{0.2 0.2}
\arr{0 2.8}{0 2.2}

\arr{-0.8 3.8}{-0.2 3.2}
\arr{-0.2 5.8}{-0.8 5.2}
\arr{0.2 5.8}{0.8 5.2}
\arr{0.8 3.8}{0.2 3.2}
\arr{-1 4.8}{-1 4.2}

\arr{1 4.8}{1 4.2}


\ellipticalarc axes ratio 1:1.8 -138 degrees from 1.3 3.9 center at 1 2.5
\arr{1.4 1.18}{1.3 1.08}

\multiput{$\alpha$} at -.7 1.6     -.7 5.6 /
\multiput{$\alpha'$} at -.7 0.4  -1.2 4.5 /
\multiput{$\beta$} at .7 1.6   2.1 2.5     .7 5.6 /
\multiput{$\beta'$} at .7 0.4  1.25 4.5 /
\multiput{$\gamma$} at -0.4 3.6  0.4 3.6 /
\multiput{$\gamma'$} at 0.2 2.7 /

\setshadegrid span <.5mm>
\hshade 3.5 -2 2  4.5 -2 2 /
\hshade 1.5 -2 2  2.5 -2 2 /
\endpicture}
$$ 
The three elements $10_1,\; 10_2\; 10_3$ are displayed in two
layers, namely in the radical layers they belong to. If we consider the
position of composition factors of the form  $10$ 
in the socle layers, we get a dual configuration, since the
subspace inside $V$ generated by the difference $10_1-10_3$
lies in the kernel of $\gamma$ and therefore belongs to $\soc_3 T(43)$.
     \medskip

Let us look at the the space 
$$
 V = 10_1\oplus 10_2\oplus 10_3,
$$
in more detail, taking into account 
all the information stored there, namely the endomorphism
$\overline\gamma = \gamma\gamma'$
as well as the images of the maps to $V$ and the kernels 
of the maps starting at $V$. One may be tempted
to look at the subspaces
$$
 \Image(\alpha'),\ \Image(\beta'),\ \Ker(\alpha),\ \Ker(\beta),
$$
however, one has to observe that the maps mentioned here are not
intrinsically given, but can be replaced by suitable others (as we have
done when we were reducing the number of parameters). For 
example, instead of looking at $\alpha'$, we have to take into account
the whole family of maps $\alpha' + c\alpha'\overline\gamma$ with $c\in k$.
Thus, the intrinsic subspaces to be considered are
$$
\begin{array}{cc}
  U_1 &= \Image(\alpha')+\Image(\alpha'\overline\gamma) = \Image(\alpha')+\Image(\overline\gamma), \cr
  U_2 &= \Image(\beta')+\Image(\overline\gamma),\qquad\qquad\qquad\qquad\ \cr
  U_3 &= \Ker(\alpha)\cap \Ker\overline\gamma, \qquad\qquad\qquad\qquad\ \cr
  U_4 &= \Ker(\beta)\cap \Ker\overline\gamma, \qquad\qquad\qquad\qquad\ 
\end{array}
$$
as well as $\Ker(\overline\gamma)$ and $\Image(\overline\gamma).$ However, 
since we see that
$$
\begin{array}{ccc}
  \Ker(\overline\gamma)& =& U_3+U_4,  \cr
  \Image(\overline\gamma) & =& U_1\cap U_2,
\end{array}
$$
it is sufficient to consider $V$ with its subspaces $U_1,\dots,U_4.$
   \medskip

This means that we deal with a vector space with four subspaces, thus with 
a representation of 
$$
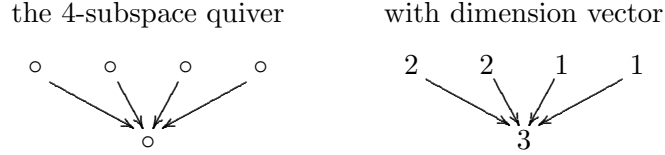
\hbox{\beginpicture
\setcoordinatesystem units <1cm,1cm>
\put{\beginpicture
\multiput{$\circ$} at 0 1  1 1  2 1  3 1  1.5 0 /
\arr{0.2 0.8}{1.3 0.2}
\arr{1.1 0.8}{1.43 0.2}
\arr{1.9 0.8}{1.57 0.2}
\arr{2.8 0.8}{1.7 0.2}
\put{the 4-subspace quiver} at 1.5 1.7
\endpicture} at 0 0
\put{\beginpicture
\put{$2$} at 0 1 
\put{$2$} at 1 1 
\put{$1$} at 2 1 
\put{$1$} at 3 1 
\put{$3$} at 1.5 0 
\arr{0.2 0.8}{1.3 0.2}
\arr{1.1 0.8}{1.43 0.2}
\arr{1.9 0.8}{1.57 0.2}
\arr{2.8 0.8}{1.7 0.2}
\put{with dimension vector} at 1.5 1.7
\endpicture} at 5 0
\endpicture}
$$
A direct calculation shows that we get the following representation:
$$
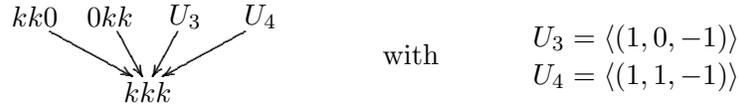
\hbox{\beginpicture
\setcoordinatesystem units <1cm,1cm>
\put{$kk0$} at 0 1 
\put{$0kk$} at 1 1 
\put{$U_3$} at 2 1 
\put{$U_4$} at 3 1 
\put{$kkk$} at 1.5 0 
\arr{0.2 0.8}{1.3 0.2}
\arr{1.1 0.8}{1.43 0.2}
\arr{1.9 0.8}{1.57 0.2}
\arr{2.8 0.8}{1.7 0.2}
\put{with} at 5 0.5
\put{$U_3 = \langle(1,0,-1)\rangle$} at 8 0.7
\put{$U_4 = \langle(1,1,-1)\rangle$} at 8 0.2
\endpicture}
$$

This is an indecomposable representation of the 4-subspace quiver, 
it belongs to a tube of rank $2$ (and is uniquely determined by its dimension vector).
Note that its endomorphism ring is a local ring of dimension 2, with radical
being the maps $V/(U_3+U_4) \to U_1\cap U_2$; and $\overline\gamma$ is just such
a map. The lattice of subspaces of $V$ generated by the subspaces $U_1,U_2,U_3,U_4$
looks as follows:
$$
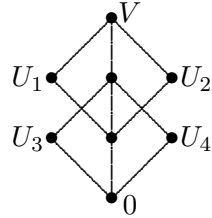
\hbox{\beginpicture
\setcoordinatesystem units <.8cm,.8cm>
\multiput{$\bullet$} at 0 0  0 1  0 2  0 3  -1 1  1 1  -1 2  1 2 /
\plot 0 0  0 3  -1 2  0 1  1 2  0 3 /
\plot 0 0  -1 1  0 2  1 1  0 0 /
\put{$0$} at 0.3 -0.1
\put{$V$} at 0.3 3.1
\put{$U_1$} at -1.4 2
\put{$U_2$} at 1.4 2
\put{$U_3$} at -1.4 1
\put{$U_4$} at 1.4 1 
\put{} at 3.3 1.5
\endpicture}
$$
Let us repeat that $\overline\gamma = \gamma\gamma'$ maps $V/(U_3+U_4)$ onto
$U_1\cap U_2$, thus we may indicate the operation of $\gamma$ and $\gamma'$ as follows:
$$
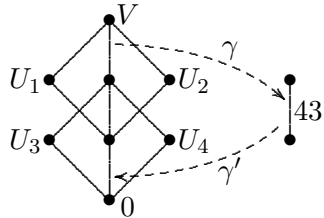
\hbox{\beginpicture
\setcoordinatesystem units <.8cm,.8cm>
\multiput{$\bullet$} at 0 0  0 1  0 2  0 3  -1 1  1 1  -1 2  1 2 /
\plot 0 0  0 3  -1 2  0 1  1 2  0 3 /
\plot 0 0  -1 1  0 2  1 1  0 0 /
\put{$0$} at 0.3 -0.1
\put{$V$} at 0.3 3.1
\put{$U_1$} at -1.4 2
\put{$U_2$} at 1.4 2
\put{$U_3$} at -1.4 1
\put{$U_4$} at 1.4 1 
\multiput{$\bullet$} at 3 1 3 2 /
\plot 3 1  3 2 /

\setdashes <1mm>
\ellipticalarc axes ratio 2:1  -62 degrees from 0.1 2.6 center at 0 1
\ellipticalarc axes ratio 2:1   62 degrees from 0.1  .4 center at 0 2
\setsolid
\arr{2.8 1.78}{2.9 1.7}
\arr{0.2 0.4}{0.1 .4}
\put{$\gamma$} at 2 2.5
\put{$\gamma'$} at 2 0.5
\put{$43$} at 3.3 1.5
\endpicture}
$$
We should stress that the last two pictures show subspace lattices (thus
composition factors are drawn as intervals between two bullets), 
in contrast to the pictures of coefficient quivers, where
the composition factors are depicted by their labels (such as
$10,\;05,\; 51,\dots$) and the lines indicate extensions of simple modules. 
	    \medskip

Note that the core of $T(43)$ is semisimple, namely of the form $10\oplus 43$,
here $10$ is just the subfactor $(U_3+U_4)/(U_1\cap U_2).$

We hope that the considerations above show well the hidden symmetries of $T(43).$

Finally, let us remark that the module $T(43)$ has a diagram $D$
in the sense of Alperin (but no strong diagram), namely the following:
$$
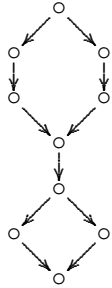
\hbox{\beginpicture
\setcoordinatesystem units <.6cm,.6cm>
\multiput{$\circ$} at 0 0  0 2  -1 4  1 4  0 6  -1 1  -1 5 
  1 1   1 5  0 3 /
\arr{-0.8 0.8}{-0.2 0.2}
\arr{-0.2 1.8}{-0.8 1.2}
\arr{0.2 1.8}{0.8 1.2}
\arr{0.8 0.8}{0.2 0.2}
\arr{0 2.8}{0 2.2}

\arr{-0.8 3.8}{-0.2 3.2}
\arr{-0.2 5.8}{-0.8 5.2}
\arr{0.2 5.8}{0.8 5.2}
\arr{0.8 3.8}{0.2 3.2}
\arr{-1 4.8}{-1 4.2}

\arr{1 4.8}{1 4.2}
\endpicture}
$$ 
(obtained from the coefficient quiver by deleting the by-path $\beta$).

\section*{10. A related algebra}

\noindent
We have used the $\Delta$-filtration of $T(43)$ in order to show that
$T(43)$ has simple socle, and this implied that the coefficient $c_{bb}$
had to be non-zero. In this way, we have obtained the somewhat strange
relations presented in the Proposition.  
Let us now consider the
same quiver $Q(10,05,51,43)$, but with the relations
$$
\begin{array}{c}
 \alpha'\alpha  = 0, \quad
 \alpha'\beta  = 0, \quad
 \beta'\alpha  = 0, \quad
 \beta'\beta = 0, \\
  \gamma'\gamma  = 0, \quad 
 \gamma'(\alpha\alpha' - \beta\beta') = 0, \quad
 (\alpha\alpha' - \beta\beta')\gamma  = 0, \quad \gamma'\alpha\alpha'\gamma = 0.
\end{array}
$$
The corresponding algebra $A'$ still is quasi-hereditary, and the
$\Delta$-modules and the $\nabla$-modules have the same shape as those
for the algebra $A = A(10,05,51,43).$ However, now it turns out that
the tilting module for $43$ is of length $11$, with a $\Delta$-filtration
of the form 
$$
 \Delta(10)\oplus\Delta(10) \mid \Delta(05)\oplus \Delta(51) \mid \Delta(43)
$$
and a similar $\nabla$-filtration.




\begin{thebibliography}{999}\frenchspacing\raggedright\small

\bibitem{alp} J.L. Alperin, Diagrams for modules, {\em J. Pure
  Appl. Algebra} \textbf{16} (1980), 111--119.

\bibitem{AK} H.H. Andersen and M. Kaneda, Rigidity of tilting modules,
  preprint.




 

\bibitem{Donkin:SA1} S. Donkin, On Schur algebras and related algebras
  I, J. Algebra 104 (1986), 310--328.

\bibitem{Donkin:SA2} S. Donkin, On Schur algebras and related algebras
  II, J. Algebra 111 (1987), 354--364.

\bibitem{Donkin:Zeit} S. Donkin, On tilting modules for algebraic
  groups, \textit{Math. Z.} \textbf{212} (1993), 39--60.

\bibitem{Donkin:book} S. Donkin, \emph{The $q$-Schur algebra},
  London Math. Soc. Lecture Note Ser., 253, Cambridge Univ. Press 1998.

\bibitem{Donkin:TiltHandbook} S. Donkin, Tilting modules for algebraic
  groups and finite dimensional algebras, \textit{Handbook of Tilting
    Theory}, 215--257, London Math. Soc. Lecture Note Ser., 332,
  Cambridge Univ. Press, Cambridge 2007.

\bibitem{Donkin:coho} S. Donkin, The cohomology of line bundles on the
  three-dimensional flag variety, \textit{J. Algebra} \textbf{307}
  (2007), 570--613.

\bibitem{D:thesis} S.R. Doty, The submodule structure of certain Weyl
  modules for groups of type $A_n$, \textit{J. Algebra} \textbf{95}
  (1985), 373--383.

\bibitem{DS} S.R. Doty and J.B. Sullivan, Filtration patterns
  for representations of algebraic groups and their Frobenius kernels,
  Math. Z. 195 (1987), 391--407.




\bibitem{DH} S.R. Doty and A.E. Henke, Decomposition of tensor
  products of modular irreducibles for $\mathrm{SL}_2$,
  \textit{Quart. J. Math.}  \textbf{56} (2005), 189--207.






\bibitem{GAP} GAP: Groups, Algorithms, Programming - a System for
Computational Discrete Algebra, \verb$http://www.gap-system.org$.

\bibitem{Green} J.A. Green, \textit{Polynomial Representations of
  $\GL_n$}, second corrected and augmented edition, Lecture Notes in
  Mathematics 830, Springer, Berlin, 2007.

\bibitem{ron} R.S. Irving, The structure of certain highest weight 
modules for SL${}_3$, \textit{J.\ Algebra} \textbf{99} (1986), 438--457.

\bibitem{Jant1} J.C. Jantzen, Darstellungen halbeinfacher 
Gruppen und kontravariante Formen, \textit{J. reine angew. Math.}
\textbf{290} (1977), 117--141.

\bibitem{Jant2} J.C. Jantzen, Weyl modules for groups of Lie type, 
\textit{Finite Simple Groups II (Proc. Durham 1978)}, Academic Press, 
London/New York (1980), 291--300.

\bibitem{Jantzen} J.C. Jantzen, \textit{Representations of Algebraic
Groups}, (2nd ed.), Mathematical Surveys and Monographs
\textbf{107}, Amer. Math. Soc., Providence 2003.

\bibitem{jgj} J. G. Jensen, On the character of some modular
  indecomposable tilting modules for $\SL_3$, \textit{J. Algebra}
  \textbf{232} (2000), 397--419.

\bibitem{kh} K. K\"{u}hne-Hausmann, Zur Untermodulstruktur der
Weylmoduln f\"{u}r $\SL_3$, \textit{Bonner Math. Schriften}, \textbf{162},
Univ. Bonn, Mathematisches Institut, Bonn 1985.

\bibitem{Martin} S.~Martin, {\em Schur Algebras and Representation
  Theory}, Cambridge Tracts in Math., vol. 112, Cambridge Univ. Press,
  Cambridge 1993.

\bibitem{Pillen} C.~Pillen, Tensor products of modules with restricted
  highest weight, \textit{Commun. Algebra} \textbf{21} (1993),
  3647--3661.


\bibitem{xi} N. Xi, Maximal and primitive elements in Weyl modules for 
type $A_2$, \textit{J.\ Algebra} \textbf{215} (1999), 735--756.

\end{thebibliography}

\begin{thebibliography}{100}\frenchspacing\raggedright\small

\bibitem[ARS]{ARS} M.Auslander, I.Reiten, S.Smal\o{}: Representation Theory
 of Artin Algebras. Cambridge Studies in Advanced Mathematics 36,
 Cambridge University Press (1995).

\bibitem[ASS]{ASS} I.Assem, D.Simson, A.Skowro\'nski: Elements of the 
 Representation Theory of Associative Algebras. 1. Techniques of
 Representation Theory. 
 London Mathematical Society Student Texts 65,
 Cambridge University Press (2007)

\bibitem[BDM]{DM} C.Bowman, S.R.Doty, S.Martin: Decomposition of tensor
  products of modular irreducible representations for $\SL_3$.

\bibitem[BR]{BR} R.A.Brualdi, H.J.Ryser: Combinatorial matrix theory.
   Encyclopedia of mathematics and its applications 39. 
   Cambridge University Press. 1991.

\bibitem[DD]{DD} M.De Visscher, S.Donkin: On projective and injective
  polynomial modules.  Math. Z. 251 (2005), 333-358.


\bibitem[DR]{DR} V.Dlab, C.M.Ringel: The module theoretic approach to
 quasi-hereditary algebras. In: Representations of Algebras
 and Related Topics (ed. H.Tachikawa, S. Brenner). London Mathematical
 Society Lectute Note Series 168. Cambridge University Press (1992).

\bibitem[M]{M} R.Mart\'inez-Villa: Algebras stably equivalent to l-hereditary.
 In: Representation Theory II (ed. V.Dlab, P.Gabriel), Springer
 Lecture Notes in Mathematics 832 (1980), 396-431.

\bibitem[R1]{R1} C.M.Ringel: Tame Algebras and Integral Quadratic Forms.
 Springer Lecture Notes in Mathematics 1099 (1984)

\bibitem[R2]{R2} C.M.Ringel: The category of modules with good filtrations 
 over a quasi-hereditary algebra has almost split sequences. 
 Math. Z. 208 (1991), 209-223. 

\bibitem[R3]{R3} C.M.Ringel: Exceptional modules are tree modules. 
   Lin. Alg. Appl. 275-276 (1998) 471-493.

\bibitem[SS]{SS}
I.Assem, D.Simson, A.Skowro\'nski: Elements of the 
 Representation Theory of Associative Algebras. 2. Tubes and Concealed
 Algebras of Euclidean Type. London Mathematical Society Student Texts 71,
 Cambridge University Press (2007)
\end{thebibliography}
\end{document}